\documentclass[12pt,a4paper,twoside,reqno]{amsart}
\usepackage{amsmath,amssymb,amsfonts,amsthm,mathrsfs}
\usepackage{enumitem}
\usepackage{times,hyperref,color}
\usepackage{graphicx}
\usepackage{cite}
\usepackage[toc,page]{appendix}
\usepackage{bm}
\usepackage{tikz}

\newcommand{\G}{\mathcal{G}}

\newcommand{\n}{{}^{(n)}}

\newcommand{\B}{\mathcal{B}}

\newcommand{\hB}{\widehat{\B}}

\newcommand{\bp}{{\bm{p}}}

\newcommand{\circfour}{\footnotesize{\mbox{\textcircled{\tiny 4}}}}
\newcommand{\circfive}{\footnotesize{\mbox{\textcircled{\tiny 5}}}}

\newcommand{\omu}{\overline{\mu}}

\newcommand{\ow}{\overline{w}}

\newcommand{\oxi}{\overline{\xi}}
\newcommand{\pp}[2]{\frac{\partial #1}{\partial #2}}
\newcommand{\vv}{\mathbf{v}}
\let\mathcal\mathscr

\newtheorem{Theorem}{Theorem}[section]

\newtheorem{Proposition}[Theorem]{Proposition}

\newtheorem{Lemma}[Theorem]{Lemma}

\theoremstyle{definition}
\newtheorem{Definition}[Theorem]{Definition}

\newtheorem{Remark}[Theorem]{Remark}

\newcommand{\oz}{\bar{z}}
\newcommand{\oZ}{\overline{Z}}

\def\hexnumber#1{\ifcase#1 0\or1\or2\or3\or4\or5\or6\or7\or8\or9\or
 A\or B\or C\or D\or E\or F\fi}

\edef\msbhx{\hexnumber\symAMSb}   % for use with amsfonts.sty

\mathchardef\emptyset="0\msbhx3F

\makeatletter
\@namedef{subjclassname@2020}{%
  \textup{2020} Mathematics Subject Classification}
\makeatother

\subjclass[2020]{32V40, 58K50, 22F50, 53A55.}

\begin{document}

%\large

\title{
Complete normal forms for real hypersurfaces in $\mathbb C^3$
\\
at $2$-nondegenerate points of Levi non-uniform rank zero
}

\author{Masoud Sabzevari}
\address{Department of Mathematics,
Shahrekord University, 88186-34141, Shahrekord, IRAN and School of
Mathematics, Institute for Research in Fundamental Sciences (IPM), 19395-5746, Tehran, IRAN}
\email{sabzevari@ipm.ir}

\date{\number\year-\number\month-\number\day}

\begin{abstract}
We construct {\it complete} normal forms for $5$-dimensional real hypersurfaces in $\mathbb C^3$ which are $2$-nondegenerate and also of Levi non-uniform rank zero at the origin point $\bp =0$. The latter condition means that the rank of the Levi form vanishes at $\bp$ but not identically in a neighborhood of it. The mentioned hypersurfaces are the only finitely nondegenerate real hypersurfaces in $\mathbb C^3$ for which their complete normal forms were absent in the literature. As a byproduct, we also treat the underlying biholomorphic equivalence problem between the hypersurfaces. Our primary approach in constructing the desired complete normal forms is to utilize the techniques derived in the theory of equivariant moving frames. It notably offers the advantage of systematic and symbolic manipulation of the associated computations.
\end{abstract}

\maketitle

\pagestyle{headings} \markright{Complete normal forms for $2$-nondegenerate real hypersurfaces in $\mathbb C^3$}

\numberwithin{equation}{section}

\section{Introduction}
\label{sec-introd}

Since 1907 that Poincar\'{e} originated in \cite{Poincare-1907} the theory of Cauchy-Riemann (CR for short) manifolds, two seminal works  \cite{Cartan-1932} and \cite{Chern-Moser} have played a profound role in the field. In the former article, Cartan laid the foundations of the systematic study of {\it equivalence problems} in CR geometry by solving it for nondegenerate hypersurfaces in $\mathbb C^2$. In the latter work \cite{Chern-Moser}, Chern and Moser introduced {\it normal forms} in CR geometry (see \cite{KKZ-17} for a survey) by constructing them for nondegenerate real hypersurfaces in arbitrary complex spaces. Although the Chern-Moser approach of normal forms mostly results in solving the underlying equivalence problems but, in comparison to Cartan's approach, it runs along a quite distinct way.

For a given real manifold $M$, let $T^cM$ be an even dimensional sub-distribution of its tangent distribution $TM$, equipped with a fiber preserving complex structure $J:T^cM\rightarrow T^cM$ with $J\circ J=-{\rm id}$. By definition \cite{BER}, $M$ is called a (abstract) {\it CR manifold} with the {\it CR structure} $T^cM$ if: a) the intersection $T^{1,0}M\cap T^{0,1}M$ is trivial where
\[
T^{1,0}M:=\big\{X-{\rm i}\,J(X): \ \ X\in T^cM \big\} \qquad {\rm and} \qquad T^{0,1}M:=\overline{T^{1,0}M}
\]
are the {\it holomorphic} and {\it anti-holomorphic} subbundles of the complexified bundle $\mathbb C\otimes TM$;
b) the holomorphic distribution $T^{1,0}M$ enjoys the {\it Frobenius condition} i.e. $[T^{1,0}M, T^{1,0}M]\subset T^{1,0}M$.

One of the most significant tools in the study of CR manifolds is the {\it Levi form} $\bf L$ which, at each point $p\in M$, is the Hermitian form ${\bf L}_p: T^{1,0}_pM\times T^{1,0}_pM\rightarrow \mathbb C\otimes\frac{T_pM}{T^c_pM}$ defined as
\[
{\bf L}_p(L_1(p), L_2(p)):={\rm i}\,[L_1,\overline L_2](p) \qquad {\rm mod} \ \mathbb C\otimes T^c_pM,
\]
for each two local vector fields $L_1$ and $L_2$ of $T^{1,0}M$ defined near $p$.
The CR manifold $M$ is called {\it Levi nondegenerate} at $p$ whenever ${\bf L}_p$ is nondegenerate as a Hermitian form.

For the significant class of real hypersurfaces of complex spaces, two subjects of equivalence problems and normal forms are widely studied and almost well-known in the nondegenerate case, specifically in the seminal works of Cartan, Chern-Moser and Tanaka \cite{Cartan-1932, Chern-Moser, Tanaka-76}.
But, in contrast, it is very little known about them around the degenerate points. For brevity and in connection with the main objective of this paper, we restrict our attention to $5$-dimensional real-analytic real hypersurfaces $M^5$ in the complex space $\mathbb C^3$ which pass through the origin point $\bp=0$\,\,---\,\,for the general case, we recommend the reader to consult \cite{BER}. In local coordinates $z_1, z_2, w:=u+{\rm i}v$ of $\mathbb C^3$, we assume that $M^5$ is represented as the graph of the real-analytic {\it real} defining function (throughout this paper, we use the notations $z$ and $\oz$ for the arrays $z_1, z_2$ and $\oz_1, \oz_2$)
\begin{equation}\label{def-fun-real}
v=v(z, \oz, u).
\end{equation}
By \cite[Theorem 4.2.6]{BER}, we can assume further that \eqref{def-fun-real} is expressed in {\it normal} coordinates, i.e. $v(z,0,u)=v(0,\oz, u)=0$.

For such a $5$-dimensional hypersurface $M^5$ and by performing necessary computations (cf. \cite[pp. 76--80]{MPS-13}), one finds the Hermitian {\it Levi matrix} associated to the Levi form $\bf L$ as
\begin{equation}\label{Lev-Matrix}
\aligned
&
{\sf L}
=
\left(\!
\begin{array}{cc}
{\rm i}
\big(
\overline{A_1}_{z_1}
+
A_1\,\overline{A_1}_u
-
{A_1}_{\oz_1}
-
\overline{A_1}\,
{A_1}_u
\big)
&
{\rm i}
\big(
\overline{A_1}_{z_2}
+
A_2\,\overline{A_1}_u
-
{A_2}_{\oz_1}
-
\overline{A_1}\,
{A_2}_u
\big)
\\
{\rm i}
\big(
\overline{A_2}_{z_1}
+
A_1\,\overline{A_2}_u
-
{A_1}_{\oz_2}
-
\overline{A_2}\,
{A_1}_u
\big)
&
{\rm i}
\big(
\overline{A_2}_{z_2}
+
A_2\,\overline{A_2}_u
-
{A_2}_{\oz_2}
-
\overline{A_2}\,
{A_2}_u
\big)
\end{array}
\!\right),
\endaligned
\end{equation}
where the functions $A_1$ and $A_2$ are defined as follows in terms of the defining function $v(z,\oz,u)$ of $M^5$
\[
\aligned
A_1
:=
-\,
\frac{v_{z_1}}{{\rm i}+v_u}, \qquad
A_2
:=
-\,
\frac{v_{z_2}}{{\rm i}+v_u}.
\endaligned
\]
Notice that at the origin point $\bp=0$, this matrix is simply
\begin{equation}\label{Levi-matrix-p}
\aligned
{\sf L}(\bp)=\left(
  \begin{array}{cc}
    v_{z_1\oz_1}(\bp) & v_{z_1\oz_2}(\bp) \\
    v_{z_2\oz_1}(\bp) & v_{z_2\oz_2}(\bp) \\
  \end{array}
\right).
\endaligned
\end{equation}
If $\sf L$ is of the full rank two at $\bp$ (and hence in a neighborhood of it), then $M^5$ is Levi nondegenerate. In the most opposite case, when this rank vanishes identically in a local neighborhood of $\bp$, then it is well-known that (\cite{BER}) $M^5$ is nothing but the {\it Levi flat} hypersurface $\mathbb C^2\times\mathbb R$ defined simply by $v=0$. Of particular interest is the intermediate scenario, when the rank of the Levi matrix $\sf L$ is not full at $\bp$ and also does not vanish identically in any local neighborhood of it.

To identify the right class of degenerate hypersurfaces in $\mathbb C^3$ with essentially unique normal form transformations (see \cite[p. 316]{Ebenfelt-98} for more details), we may introduce the extra assumption of $2$-nondegeneracy defined as below. Let
\[
Q(z,\oz, w,\ow)=0
\]
be the complex defining equation of the hypersurface $M^5$ which is plainly obtained by setting $u=\frac{w+\ow}{2}$ and $v=\frac{w-\ow}{2\rm i}$ in its real counterpart \eqref{def-fun-real}. Let $\nabla^hQ$ denotes the holomorphic vector gradients of $Q$. Then, by definition \cite{BER, Ebenfelt-98}, the hypersurface $M^5$ is $2$-{\it nondegenerate} at $\bp$ if the set of triple vectors
\[
\big\{\nabla^hQ(\bp, \overline\bp), \qquad \overline{\mathcal L}_j\big(\nabla^hQ(\bp, \overline\bp)\big), \qquad \overline{\mathcal L}_j\overline{\mathcal L}_k\big(\nabla^hQ(\bp, \overline\bp)\big), \qquad j,k=1,2\big\}
\]
spans $\mathbb C^3$, on the contrary of its proper subset $\{\nabla^hQ(\bp, \overline\bp), \mathcal L_j\big(\nabla^hQ(\bp, \overline\bp)\big), \ j=1,2\}$\,\,---\,\,if the latter vectors span $\mathbb C^3$, then $M^5$ is nondegenerate at $\bp$. Here, $\overline{\mathcal L}_1$ and $\overline{\mathcal L}_2$ are the generators of the anti-holomorphic distribution $ T^{0,1}M^5$, i.e.
\[
\overline{\mathcal L}_j=\frac{\partial}{\partial\oz_j}+\frac{\partial\overline Q}{\partial\oz_j}\,\frac{\partial}{\partial\ow}, \qquad j=1,2.
\]

Over the past few decades, there has been significant research devoted to the study of $5$-dimensional $2$-nondegenerate real hypersurfaces in $\mathbb{C}^3$. Notably, the work \cite{Ebenfelt-98} of Ebenfelt stands out, where he made substantial contributions by constructing normal forms for these hypersurfaces within two distinct subclasses. The first subclass comprises hypersurfaces whose associated Levi matrix ${\sf L}(\bp)$ is of rank one\,\,---\,\,or equivalently ${\sf L}(\bp)$ admits one zero and one nonzero eigenvalue.  Ebenfelt showed in the first part of \cite[Theorem {\bf A}]{Ebenfelt-98} that after applying suitable normalizations, each hypersurface of this kind can be transformed into one of the following {\it partial} normal forms
\begin{equation}\label{Ebenfelt-partial-NF-i}
\aligned
{\rm (A.i.1)}&: \qquad v=z_1\oz_1+z_2^2\oz_2+z_2\oz_2^2+\gamma\,\big(z_1^2\oz_2+z_2\oz_1^2\big)+{\rm O}(|z|^4)+{\rm O}(|u||z|^2), \qquad \gamma=0,1,
\\
{\rm (A.i.2)}&: \qquad v=z_1\oz_1+z_1^2\oz_2+z_2\oz_1^2+{\rm O}(|z|^4)+{\rm O}(|u||z|^2),
\\
{\rm (A.i.3)}&: \qquad v=z_1\oz_1+z_1z_2\oz_2+z_2\oz_1\oz_2+{\rm O}(|z|^4)+{\rm O}(|u||z|^2).
\endaligned
\end{equation}
He also showed that (\cite[Theorem 4.2.8]{Ebenfelt-98}) $M^5$ is biholomorphically equivalent to a normal form hypersurface of the type (A.i.2) if and only if it is {\it Levi uniform of rank} $1$, meaning that the rank of the Levi matrix $\sf L$ is constantly equal to $1$ in a neighborhood of $\bp$. The class of such real hypersurfaces, sometimes denoted by $\frak C_{2,1}$, is studied extensively among the literature. For instance, Ebenfelt investigated in \cite{Ebenfelt-01} their Cartan geometry and Isaev-Zaitsev, Medori-Spiro and Merker-Pocchiola studied in \cite{Isaev-Zaitsev-13, Medori-Spiro-14, Merker-Pocchiola-20} their biholomorphic equivalence problem through Cartan's classical approach. Moreover, Fels and Kaup in \cite{Fels-Kaup-08} classified homogeneous hypersurfaces belonging to this class.

The second subclass of $2$-nondegenerate hypersurfaces studied in \cite{Ebenfelt-98} comprises real hypersurfaces in $\mathbb C^3$ whose  associated Levi matrix $\sf L$ vanishes at $\bp$ but not in a neighborhood thereof. In this scenario, we refer the $2$-nondegenerate point $\bp$ as having the {\it Levi non-uniform} rank zero (cf. \cite{Ebenfelt-01}). In the second part of Theorem {\bf A} of \cite{Ebenfelt-98}, Ebenfelt showed that such hypersurfaces can be brought to the following biholomorphically inequivalent types of {\it partial} normal forms
\begin{equation}\label{Ebenfelt-Branches}
\aligned
{\rm (A.ii.1)} \qquad v&=z_1z_2\oz_1+z_1\oz_1\oz_2+r\,\big(z_1^2\oz_2+z_2\oz_1^2\big)+{\rm O}(|z|^4)+{\rm O}(|u||z|^2),
\\
{\rm (A.ii.2)} \qquad v&=z_1z_2\oz_1+z_1\oz_1\oz_2+z_1^2\oz_2+z_2\oz_1^2+{\rm i}\,\big(z_1^2\oz_1-z_1\oz_1^2\big)+{\rm O}(|z|^4)+{\rm O}(|u||z|^2),
\\
{\rm (A.ii.3)} \qquad v&=z_1z_2\oz_1+z_1\oz_1\oz_2+z_2^2\oz_1+z_1\oz_2^2+\lambda\,z_2^2\oz_2+\overline\lambda\, z_2\oz_2^2+{\rm O}(|z|^4)+{\rm O}(|u||z|^2),
\\
{\rm (A.ii.4)} \qquad v&=z_1^2\oz_1+z_1\oz_1^2+z_2^2\oz_2+z_2\oz_2^2+\sigma\,z_1^2\oz_2+\overline\sigma\,z_2\oz_1^2+\nu\,z_2^2\oz_1+\overline\nu\,z_1\oz_2^2
\\
&\ \ \ \ \ \ \ \ \ \ \ \ +{\rm O}(|z|^4)+{\rm O}(|u||z|^2),
\\
{\rm (A.ii.5)} \qquad v&=\eta\,z_1^2\oz_1+\overline\eta z_1\oz_1^2+z_1^2\oz_2+z_2\oz_1^2+z_2^2\oz_1+z_1\oz_2^2+{\rm O}(|z|^4)+{\rm O}(|u||z|^2),
\endaligned
\end{equation}
with $r\in\mathbb R$ and $\lambda,\sigma,\nu,\eta\in\mathbb C$ where $r>0$, $\lambda\neq 0$ and $\sigma\nu\neq 1$.

By definition (\cite{Ebenfelt-98}), a normal form associated with a hypersurface $M^5$ is {\it complete}\,\,---\,\,in contrast to {\it partial}\,\,---\,\,if the corresponding transformation of $M^5$ to the normal form is {\it unique} modulo a finite dimensional choice of normalizations. Neither of the groups of normal forms \eqref{Ebenfelt-partial-NF-i} and \eqref{Ebenfelt-Branches} are complete. However, Ebenfelt in \cite[Theorem {\bf B}]{Ebenfelt-98} developed the former partial forms \eqref{Ebenfelt-partial-NF-i} into their complete counterparts by applying infinitely many appropriate normalizations. In addition, recently in \cite{FMT-22, KK-22}, two {\it convergent} complete normal forms are constructed in the case (A.i.2).

But, completing the second group \eqref{Ebenfelt-Branches} of Ebenfelt's partial normal forms has been overlooked among the literature.
From a computational perspective, the removal of second order monomials $z_j\oz_k$, $j,k=1,2$ in the expressions \eqref{Ebenfelt-Branches} makes the construction of their corresponding (complete) normal forms more challenging.

The main objective of this paper is to develop the partial normal forms \eqref{Ebenfelt-Branches} to their complete forms. Let us outline the extent to which we aim to apply the requisite normalizations in order to attain the desired {\it uniqueness} in the above mentioned definition of a complete normal form. In \cite{Ershova-01}, Ershova computed the sharp upper bounds for dimensions of the isotropy groups associated to the hypersurfaces of types (A.i.1)--(A.i.3) and (A.ii.1)--(A.ii.5) at $\bp$. She showed that the maximum dimensions are enjoyed by the so-called  {\it model hypersurfaces} associated to each type. Our strategy in this paper is to apply the requisite normalizations until we succeed in reducing the number of remaining unnormalized group parameters to a finite value which is less than or equal to Ershova's upper bound.

As suggested in \cite{Olver-2018}, our primary approach toward constructing the desired complete normal forms is to utilize the techniques derived in the theory of {\it equivariant moving frames} (see $\S$\ref{sec-prel-norm} below). This theory, which is initiated and developed by Peter Olver and his school (\cite{Olver-Fels-99, Olver-Pohjanpelto-05, Olver-Pohjanpelto-08, Olver-2015}) in the recent three decades, is indeed a modern and far-reaching reformulation of Cartan's classical approach to moving frames. Besides its various applications (see \cite{Olver-2015} and the reference therein), it notably provides a {\it systematic} and {\it algorithmic} way, based on {\it symbolic computations}, to {\it simultaneously} construct the normal forms (\cite{Olver-2018}) and solve their underlying equivalence problems (\cite{Valiquette-SIGMA}). It provides a concrete and striking bridge between Cartan's classical approach in solving equivalence problems and the theory of normal forms. It may answer in part the two questions $Q^{\circfour}$ and $Q^{\circfive}$ introduced in \cite[p. 257]{FMT-22}.

First applications of the equivariant moving frames in CR geometry appeared recently in \cite{OSV-23, Normal-Form}, where the normal forms of $3$-dimensional {\it nondegenerate} real hypersurfaces in $\mathbb C^2$ and $5$-dimensional {\it totally nondegenerate} CR surfaces in $\mathbb C^4$ are constructed. These works highlight the emerging role of equivariant moving frames in CR geometry and exhibit their potential for addressing various problems in this field. The current work is the first application of the equivariant moving frames in the degenerate case.

This paper is organized as follows. In the next Section \ref{sec-prel-norm}, we prepare requisite materials to launch the equivariant moving frame method for constructing the desired normal forms. In Section \ref{sec-elementar-nor}, we apply partial normalizations in the lower orders $\leq 3$. One observes that while the rank zero Levi non-uniformity at $\bp$ prohibits any normalization in order two, normalizations in order three depends upon certain circumstances imposed by the $2$-nondegeneracy condition. This forces us to apply the next normalizations along the five distinct branches \eqref{Ebenfelt-Branches}. In the next three sections \ref{sec-A.ii.1-A.ii.2}, \ref{sec-A.ii.3} and \ref{sec-A.ii.4-A.ii.5}, we endeavor to complete these five partial normal forms. For this purpose, we need to pursue in each branch the final set of normalizations in order six.  As a byproduct, we independently find Ershova's model hypersurfaces \cite{Ershova-01} as those admitting the minimum number of normalizations or equivalently as those having the maximum dimension of isotropy groups at $\bp$ (cf. \cite{Valiquette-SIGMA}). Finally, in the short section \ref{sec-equiv-prob}, we investigate the biholomorphic equivalence problem between $2$-nondegenerate hypersurfaces, considered in this paper.

\section{Preliminary materials}
\label{sec-prel-norm}

We aim in this section to prepare requisite materials toward launching the equivariant moving frame method for constructing the desired normal forms. It necessitates first to consider the underlying pseudo-group action of holomorphic transformations.

\subsection{Holomorphic pseudo-group}

In local coordinates $z_1, z_2, w$, the pseudo-group $\mathcal G$ of origin-preserving holomorphic transformations of $\mathbb C^3$ consists of diffeomorphisms
\[
(z_1, z_2, w)\rightarrow (Z_1(z_1, z_2,w), Z_2(z_1, z_2,w), W(z_1, z_2,w))
\]
which enjoy the Cauchy-Riemann equations
\[
\frac{\partial Z_j}{\partial \overline z_k}=\frac{\partial Z_j}{\partial \overline w}=\frac{\partial W}{\partial \overline z_k}=\frac{\partial W}{\partial \overline w}=0, \qquad j,k=1,2.
\]

Expanding $w=u+{\rm i} v$ and accordingly $W(z,w):=U(z,\oz,u,v)+{\rm i} V(z,\oz,u,v)$ to the real and imaginary parts, one finds equivalently the first order {\it determining equations} of $\mathcal G$ as (cf. \cite{Normal-Form})
\begin{equation}
\label{eq: determining equations}
\aligned
&\frac{\partial Z_j}{\partial \overline z_k} = \frac{\partial \oZ_j}{\partial z_k}= 0,\qquad \frac{\partial Z_j}{\partial v}={\rm i} \frac{\partial Z_j}{\partial u}, \qquad \frac{\partial \oZ_j}{\partial v}=-{\rm i} \frac{\partial \oZ_j}{\partial u},
\\
&\frac{\partial U}{\partial z_j} = {\rm i}\frac{\partial V}{\partial z_j},\qquad \frac{\partial U}{\partial \oz_j} = -{\rm i}\frac{\partial V}{\partial \oz_j}, \qquad \frac{\partial U}{\partial v}=-\frac{\partial V}{\partial u}, \qquad \frac{\partial V}{\partial v}=\frac{\partial U}{\partial u}, \qquad j, k= 1, 2.
\endaligned
\end{equation}

It follows that the infinitesimal counterpart of $\mathcal G$, namely the Lie algebra $\frak g:=\frak{hol}(\mathbb C^3, 0)$ of local holomorphic automorphisms of $\mathbb C^3$ around the origin, consists of the real vector fields
\begin{equation}\label{eq: v}
\vv = \sum_{j=1}^2\,\xi^j(z,u,v)\pp{}{z_j} +\sum_{j=1}^2\, \overline{\xi}^j(\oz,u,v)\pp{}{\oz_j} + \eta(z,\oz,u,v)\pp{}{u} + \phi(z,\oz,u,v)\pp{}{v},
\end{equation}
with $\overline{\xi}^j(\oz,u,v):=\overline{\xi^j(z,u,v)}$, which enjoys the second order {\it infinitesimal determining equations}
\begin{equation}\label{eq: infinitesimal determining equations}
\begin{gathered}
\xi^j_{\oz_k}=\oxi^j_{z_k}=0,\qquad \xi^j_{v} = {\rm i}\,\xi^j_{u},\qquad \oxi^j_{v}=-{\rm i}\,\oxi^j_{u},
 \\
\phi_{z_k} = -{\rm i}\,\eta_{z_k},\qquad \phi_{\oz_k} = {\rm i}\,\eta_{\oz_k},\qquad \phi_{u} = - \eta_{v},\qquad \phi_{v}=\eta_{u},
\\
\xi^j_{\oz_k,a} = 0,\qquad \xi^j_{z_kv} = {\rm i}\,\xi^j_{z_ku},\qquad \xi^j_{uv} = {\rm i}\,\xi^j_{uu},\qquad \xi^j_{vv} = -\xi^j_{uu},
\\
\oxi^j_{z_k,a} = 0, \qquad  \oxi^j_{\oz_k v} = -{\rm i}\,\oxi^j_{\oz_k u},\qquad \oxi^j_{uv} = -{\rm i}\,\oxi^j_{uu},\qquad \oxi^j_{vv} = -\oxi^j_{uu},
\\
\eta_{z_j\oz_k}=0,  \qquad \eta_{vv}=-\eta_{uu}, \qquad \eta_{z_kv}={\rm i}\, \eta_{z_ku},\qquad
\eta_{\oz_k v} = -{\rm i}\, \eta_{\oz_k u},
\\
\phi_{z_k,a} = -{\rm i}\,\eta_{z_k,a},\qquad \phi_{\oz_k,a} = {\rm i}\, \eta_{\oz_k,a},\qquad \phi_{u,a} = -\eta_{v,a},\qquad
\phi_{v,a} = \eta_{u,a}.
\end{gathered}
\end{equation}
Higher order infinitesimal determining equations are obtained by applying further derivations on these equations.

One may consider the induced action of $\mathcal G$ on $\mathbb C^3$ on a given $5$-dimensional real hypersurface $M^5\subset\mathbb C^3$, represented in local coordinates as the graph of some defining function $v=v(z, \oz, u)$. For $0\leq n\leq\infty$, let ${\rm J}^n:={\rm J}^n(M^5)$ denotes the associated $n$-th order jet space of $M^5$. In local coordinates, it comprises the tuples $v^{(n)}=(z,\oz,u,v,\ldots, v_J, \ldots)$, where $J$ is a symmetric multi-index of order $\leq n$ with entries $j_k=z_1, z_2, \oz_1, \oz_2, u$. As explained in \cite{Olver-1995}, the action of $\mathcal G$ on $M^5$ can be extended to action of the {\it prolonged} group $\mathcal G^{(n)}$ on the jet space ${\rm J}^n$. In local coordinates, $\mathcal G^{(n)}$ consists of the derivatives of the transformations in $\mathcal G$, up to order $n$. From now on and for $j=1,2$, we denote by $Z_j, \oZ_j, U, V$ the {\it lift} of the source variables $z_j, \oz_j, u, v$ under prolonged transformations. In a similar fashion, transformation of the above jet point $v^{(n)}$ is denoted by $V^{(n)}=(Z,\oZ,U,V,\ldots, V_J, \ldots)$, where here the entries of the multi-indices $J$ are $Z_1, Z_2, \oZ_1, \oZ_2, U$. In \cite{Olver-1995}, the formula for explicitly computing the target variables $V_J$ can be found. However, in this paper, we do not need to compute them explicitly as they will be considered from a {\it symbolic} standpoint.

Infinitesimally, the already mentioned prolonged action $\mathcal G^{(\infty)}$ on $J^{\infty}$ is determined by prolonging the above vector field $\bf v$, namely by
\begin{equation}
\label{v-infty}
\vv^{(\infty)} = \sum_{j=1}^2 \xi^j\, \frac{\partial}{\partial z_j} + \sum_{j=1}^2\oxi^j \frac{\partial}{\partial \oz_j} + \eta \frac{\partial}{\partial u} + \sum_{\sharp J\geq 0} \phi^{J}\frac{\partial}{\partial v_J},
\end{equation}
where the vector components $\phi^{J}$\,\,---\,\,not to be confused with the derivations $\phi_J$\,\,---\,\,are defined recursively by the prolongation formula
\begin{equation}\label{eq: prolongation formula 2}
\phi^{J, a}=D_{ a}\phi^{J}-(D_{ a}\xi^j) \,v_{z_j, J}-(D_{a}\overline\xi^j) \,v_{J, \oz_j}-(D_{a}\eta)\, v_{J, u_j},
\end{equation}
for $a=z_1, z_2, \oz_1, \oz_2, u$. Here, $D_a$ is the total differentiation operator with respect to $a$ (cf. \cite{Olver-1995}).

To the vector components of $\bf v$, we associate the zeroth order Maurer-Cartan forms $\mu^1, \mu^2, \alpha, \gamma$ with the assignments
\begin{equation}\label{MC-vector}
\mu^j\leftrightarrow\xi^j, \qquad \omu^j\leftrightarrow\oxi^j, \qquad \alpha\leftrightarrow\eta, \qquad \gamma\leftrightarrow\phi, \qquad {\rm for} \ j=1,2.
\end{equation}
These differential forms can be extended to higher orders $\mu^j_J, \omu^j_J, \alpha_J, \gamma_J$, where the entries of the multi-indices $J$ are again $Z_1, Z_2, \oZ_1, \oZ_2, U$.  The coordinate expressions of the Maurer-Cartan forms can be found in \cite{Olver-Pohjanpelto-05}, although in the subsequent computations, we treat them symbolically and therefore do not require their explicit forms here.

The Maurer-Cartan forms $\mu^j_J, \omu^j_J, \alpha_J, \gamma_J$ with $0\leq\#J\leq\infty$ are not linearly independent.
Detecting their dependencies, one may apply Theorem 6.1 of \cite{Olver-Pohjanpelto-05} which, according to the infinitesimal determining equations \eqref{eq: infinitesimal determining equations}, gives in orders $\leq 2$ that
\begin{equation}\label{MC forms=relations}
\begin{gathered}
\mu^j_{\oZ_k}=\omu^j_{Z_k}=0,\qquad \mu^j_{V} = {\rm i}\,\mu^j_{U},\qquad \omu^j_{V}=-{\rm i}\,\omu^j_{U},
 \\
\gamma_{Z_k} = -{\rm i}\,\alpha_{Z_k},\qquad \gamma_{\oZ_k} = {\rm i}\,\alpha_{\oZ_k},\qquad \gamma_{U} = - \alpha_{V},\qquad \gamma_{V}=\alpha_{U},
\\
\mu^j_{\oZ_k,a} = 0,\qquad \mu^j_{Z_kV} = {\rm i}\,\mu^j_{Z_kU},\qquad \mu^j_{UV} = {\rm i}\,\mu^j_{UU},\qquad \mu^j_{VV} = -\mu^j_{UU},
\\
\omu^j_{Z_k,a} = 0, \qquad  \omu^j_{\oZ_k V} = -{\rm i}\,\omu^j_{\oZ_k U},\qquad \omu^j_{UV} = -{\rm i}\,\omu^j_{UU},\qquad \omu^j_{VV} = -\omu^j_{UU},
\\
\alpha_{Z_j\oZ_k}=0,  \qquad \alpha_{VV}=-\alpha_{UU}, \qquad \alpha_{Z_kV}={\rm i}\, \alpha_{Z_kU},\qquad
\alpha_{\oZ_k V} = -{\rm i}\, \alpha_{\oZ_k U},
\\
\gamma_{Z_k,a} = -{\rm i}\,\alpha_{Z_k,a},\qquad \gamma_{\oZ_k,a} = {\rm i}\, \alpha_{\oZ_k,a},\qquad \gamma_{U,a} = -\alpha_{V,a},\qquad
\gamma_{V,a} = \alpha_{U,a},
\end{gathered}
\end{equation}
for $a\in\{Z_j, \oZ_j, U, V, \ j=1,2\}$.
One finds the higher order linear relations by applying successive derivations on the both sides of these equalities.
Accordingly, we have the following basis of Maurer-Cartan forms
\begin{equation}\label{MC-1}
\mu^j_{Z^\ell U^k}, \qquad  \omu^j_{\oZ^\ell U^k}, \qquad  \alpha_{U^j V}, \qquad \alpha_{Z^\ell U^k}, \qquad \alpha_{\oZ^\ell U^k}, \qquad \gamma
\end{equation}
for $j=1,2$, $k\in\mathbb N_0$ and for $\ell:=(l_1, l_2)=\mathbb N_0^2$, where by $Z^\ell$ we mean $Z_1^{l_1} Z_2^{l_2}$ (we denote $\mathbb N_0$ the collection of natural numbers added by $0$).

\subsection{Moving frame}

For $0\leq n\leq\infty$, let $\mathcal B^{(n)}$ denotes the $n$-th order {\it lifted fiber bundle} which, in local coordinates, is parameterized by the pairs $(\varphi^{(n)}, v^{(n)})\in\mathcal G^{(n)}\times {\rm J}^n$. The action $\mathcal G$ can be naturally extended on the lifted bundle $\mathcal B^{(n)}$ by the right composition of the pseudo-group jets
\begin{equation}\label{eq: right action}
\psi\cdot (\varphi^{(n)}, v^{(n)})=((\varphi\circ\psi^{-1})^{(n)}, V^{(n)}), \qquad {\rm for \ any} \ \psi\in\mathcal G.
\end{equation}

\begin{Definition}
A \emph{partial right moving frame} of  order $n$ is a right-invariant local subbundle $\hB\n\subset \B\n$, meaning that $\psi\cdot\hB\n \subset \hB\n$ for all $\psi \in \G$ where the right action \eqref{eq: right action} is defined. If the subbundle $\hB\n$ forms the graph of a right-invariant section of $\B\n$, it defines an equivariant moving frame.
\end{Definition}

The construction of a partial moving frame heavily relies on selecting an appropriate {\it cross-section} to the corresponding action. In the most practical and standard way, the so-called {\it coordinate} cross-section is obtained by setting some appropriate $V_J$, called the {\it phantom invariants}, to some suitable constants $c_J$. Solving the resulting equations
\begin{equation}
\label{normalization-eq}
dV_J=c_J,
\end{equation}
which are known as  the {\it normalization equations}, reveals the normalized expressions of the group parameters in terms of the jet coordinates. Inserting back these expressions into the original expressions of the non-phantom functions $V_J$, one obtains a complete set of differential invariants associated with the action of $\mathcal G$ on the manifolds. This observation motivates the term {\it lifted differential invariants} for the target jet coordinates $V_J$ in the context of equivariant moving frames. (\cite{Olver-Fels-99, Olver-Pohjanpelto-08}).

Associated to the local coordinates $z_j, \oz_j, u, j=1,2$ of $M^5$, we have {\it invariant horizontal forms} $\omega^{Z_j}, \omega^{\oZ_j}=\overline{\omega^{Z_j}}$ and $\omega^U$ which actually are the {\it lifts} of the standard $1$-forms $dz_j, d\oz_j, du$ under the group action $\mathcal G$. As before, the coordinate expressions of these forms are plainly accessible (cf. \cite{Olver-Pohjanpelto-08}) but we do not need them in this paper.

In general, substituting the resulted expressions of the group parameters among the process of cross-section normalizations into the prolonged transformation formulae establishes what is known as the {\it invariantization} process which maps each differential function, differential form, differential operator and etc to its invariant counterpart (\cite{Olver-Fels-99, Olver-Pohjanpelto-08}). We denote this invariantization operator by $\iota$.

\subsection{Recurrence formula}

The most powerful tool in the theory of equivariant moving frames is the so-called recurrence formula. In general, the invariantization map $\iota$ and the exterior differentiation $d$ do not commute and this formula measures to what extent they may differ.

\begin{Theorem}
 ({\it cf.} \cite[Theorem 25]{Olver-Pohjanpelto-08}) If $\Omega$ is a differential form defined on ${\rm J}^\infty$ then it enjoys the {\sl recurrence relation}
\begin{equation*}
d\iota(\Omega)=\iota\big[d\Omega+{\bf v}^{(\infty)}(\Omega)\big]
\end{equation*}
 where ${\bf v}^{(\infty)}$ is the prolonged infinitesimal vector field \eqref{v-infty} and ${\bf v}^{(\infty)}(\Omega)$ denotes the Lie derivative of $\Omega$ along it.
\end{Theorem}

Accordingly, the recurrence relations in our case are represented as
\begin{equation}
\label{rec-formula}
\aligned
dZ_j&=\omega^{Z_j}+\mu^j, \qquad d\oZ_j=\omega^{\oZ_j}+\overline\mu^j,
\\
dU&=\omega^{U}+\alpha,
\\
dV_J&=\varpi_{J}+{\widehat{\phi}}^{J},
\endaligned
\end{equation}
for $j=1,2$ and $\sharp J\geq 0$, where from now on we denote
\[
\varpi_{J}:=V_{J, Z_j}\,\omega^{Z_j}+V_{J, \oZ_j}\,\omega^{\oZ_j}+V_{J, U}\,\omega^{U}
\]
and where ${\widehat{\phi}}^{J}$ is the invariantization of the vector coefficient $\phi^{J}$ of ${\bf v}^{\infty}$ which, roughly speaking, is obtained by substituting in the expression of $\phi^J$ the source jet variables $v_J$ by their lifted counterparts $V_J$ and the vector components $\xi^j, \eta, \phi$ with their corresponding Maurer-Cartan forms  \eqref{MC-vector}.

In practice, as the order of $J$ increases, the expressions of the lifted invariants $V_J$ can grow explosively. However, the recurrence formula offers a {\it systematic} and {\it symbolic} approach to normalize the corresponding Maurer-Cartan forms instead of directly normalizing the group parameters by solving the aforementioned normalization equations. More precisely, if $V_J=c_J$ is a normalization equation for a phantom lifted invariant $V_J$ and constant $c_J$, we apply it in the corresponding recurrence relation $dV_J$ and solve the resulting equation for some suitable Maurer-Cartan form. It notably offers the advantage of performing linear algebraic manipulations for normalizing the corresponding Maurer-Cartan forms, without the need for explicit expression of the lifted invariants, Maurer-Cartan forms, normalized expressions of the group parameters and etc.

\subsection{Normal form}
 \label{sec-normal-form}

Let $M^5$ be a real-analytic $5$-dimensional real hypersurface in $\mathbb C^3$, passing through the origin point $\bp=0$ and represented in local coordinates $z_1, z_2, w=u+{\rm i}v$ as the graph of some defining equation $v:=v(z,\oz,u)$. For a $5$-tuple $J=(j_1, j_2, k_1, k_2, l)$, denote $v_J:=v_{z_1^{j_1}z_2^{j_2}\oz_1^{k_1}\oz_2^{k_2}u^l}$,  $x^J:=z_1^{j_1}z_2^{j_2}\oz_1^{k_1}\oz_2^{k_2}u^l$ and $J!:=j_1!\, j_2!\, k_1!\, k_2!\, l!$. Then, around $\bp$, the Taylor series of $M^5$ is
\begin{equation}\label{Taylor-ser}
v(z,\oz,u)=\sum_{\sharp J\geq 0} \, \frac{v_{J}(\bp)}{J!}\,x^J.
\end{equation}
 As suggested in \cite{Olver-2018}, we identify this Taylor series expansion with the restriction $v^{(\infty)}|_{\bp}$ of the jet coordinates to the point $\bp$. A normal form of $M^5$ is made by employing the group transformations of $\mathcal G$ to simplify, as much as possible, the coefficients of the above Taylor series. This process is analogous to finding some practical partial cross-section\,\,---\,\,or equivalently a moving frame\,\,---\,\,associated with the action of the pseudo-group $\mathcal G$. Along the process of normalizations, the jet coordinates $v_J$ transform into differential invariants $V_J$ which no longer  depend on any potentially remaining unnormalized group parameters (cf. \cite{Valiquette-SIGMA}). It amounts to converting the Taylor series \eqref{Taylor-ser} to the {\it normal form}
 \begin{equation}\label{NF-general}
 v(z,\oz,u)=\sum_{\sharp J\geq 0} \, \frac{V_{J}(\bp)}{J!}\,x^J.
 \end{equation}

It follows that the essential conditions that determine the coefficients of the space of the desired normal forms are actually imposed by the corresponding cross-section through the normalization equations \eqref{normalization-eq}. We refer the reader to  \cite{Olver-2018} for further relevant information and details.

\section{Elementary normalizations}
\label{sec-elementar-nor}

Now we are ready to launch the equivariant moving frame technique in constructing complete normal forms for real hypersurfaces $M^5\subset\mathbb C^3$ at their $2$-nondegenerate point $\bp =0$ of non-uniformly Levi rank zero. Benefiting the powerful recurrence formula \eqref{rec-formula}, we manage our computations symbolically and order by order without requiring any explicit coordinate expression. Since the lifted invariants and Maurer-Cartan forms respect the conjugation relation, i.e.
\[
\overline{V_{Z^{\ell_1}\oZ^{\ell_2}U^l}}=V_{Z^{\ell_2}\oZ^{\ell_1}U^l}, \qquad {\rm and} \qquad \overline{\mu^j_{Z^{\ell_1}\oZ^{\ell_2}U^l}}=\omu^j_{Z^{\ell_2}\oZ^{\ell_1}U^l}, \qquad \overline{\alpha_{Z^{\ell}U^l}}=\alpha_{\oZ^{\ell}U^l},
\]
then all our upcoming computations and normalizations will respect this relation, as well. Accordingly, we will consider the normalization of only one of the two conjugative Maurer-Cartan forms.
For brevity, we occasionally employ Einstein's summation convention to present our expressions.

\subsection{Orders zero and one}

In light of \eqref{rec-formula}, applying the recurrence formula on the zeroth and first orders lifted invariants gives
\[
\aligned
dZ_j&=\omega^{Z_j}+\mu^j,
\\
dU&=\omega^U+\alpha, \ \ \ \, \qquad dV=\varpi+\gamma
\\
dV_{Z_j}&=\varpi_{Z_j}-\mu^k_{Z_j} V_{Z_k}-\alpha_{Z_j} V_U-{\rm i}\,\alpha_{Z_j},
\\
dV_U&=\varpi_U-\mu^k_U V_{Z_k}-\omu^k_U V_{\oZ_k}-\alpha_U V_U-\alpha_V
\endaligned
\]
for $j=1,2$. Thus, by selecting the orders zero and one cross-section $Z_j=U=V=V_{Z_j}=V_U=0$ (with the same conjugations), the above relations can be solved for the Maurer-Cartan forms
\[
\mu^j=-\omega^{Z_j}, \qquad \alpha=-\omega^U, \qquad \gamma=0, \qquad \alpha_{Z_j}=-{\rm i}\,\varpi_{Z_j}, \qquad    \alpha_V=\varpi_U.
\]

Analysing carefully the prolongation formula \eqref{eq: prolongation formula 2} leads us to the following general observation.

\begin{Lemma}
\label{lem-ord-one}
Let $j\geq 0$ and $(0,0)\neq\ell\in\mathbb N_0^2$. Then
\begin{itemize}
  \item[$1)$] it is possible to normalize the Maurer-Cartan form $\alpha_{Z^\ell U^j}$ by setting $V_{Z^\ell U^j}\equiv0$.
  \item[$2)$] it is possible to normalize the Maurer-Cartan form $\alpha_{U^jV}$ by setting $V_{U^{j+1}}\equiv0$.
\end{itemize}
\end{Lemma}

\proof
We establish the proof for the first assertion, noting that the proof for the second assertion follows along a similar argument. In order to consider the recurrence relation of $V_{Z^\ell U^l}$, we shall inspect the vector component $\phi^{z^\ell u^l}$ in the prolonged vector field \eqref{v-infty}. As $\ell=(l_1, l_2)\neq (0,0)$, then without loss of generality we assume that $l_1\neq 0$. One plainly sees that
\[
\phi^{z_1}=\phi_{z_1}+\cdots
\]
where $"\cdots"$ stands for some combinations of the first order jets $v_{a}$, $a=z_1, z_2, \oz_1, \oz_2, u$ and the derivations of the other vector components. Then, thanks to the prolongations formula \eqref{eq: prolongation formula 2}, one readily verifies for every ${\sf a}=z_1, z_2, u$ that
\[
\aligned
\phi^{z_1 {\sf a}}&=D_{\sf a}\,\phi^{z_1}-(D_{\sf a}\xi^j)\,v_{z_1z_j}-(D_{\sf a}\oxi^j)\,v_{z_1\oz_j}-(D_{\sf a}\eta)\,v_{z_1u}
\\
&=\phi_{z_1\sf a}+\phi_{z_1v}\,v_{\sf a}-(D_{\sf a}\xi^j)\,v_{z_1z_j}-(D_{\sf a}\oxi^j)\,v_{z_1\oz_j}-(D_{\sf a}\eta)\,v_{z_1u}+\cdots.
\endaligned
\]
By \eqref{eq: infinitesimal determining equations}, we have $\phi_{z_1\sf a}=-{\rm i}\eta_{z_1\sf a}$, which after invariantization yields the following  recurrence relation
\[
dV_{Z_1\sf A}=-{\rm i}\,\alpha_{Z_1\sf A}+\cdots
\]
where $\sf A$ is the lift of the source variable $\sf a$ and  where $\alpha_{Z_1\sf A}$ does not appear in the $"\cdots"$ part. Now, applying a simple induction shows that the recurrence relation of every arbitrary $V_{Z_1^{j+1} Z_2^k U^l}$ is of the form
\[
dV_{Z_1^{j+1} Z_2^k U^l}=-{\rm i}\, \alpha_{Z_1^{j+1} Z_2^k U^l}+\cdots
\]
where the Maurer-Cartan form $\alpha_{Z_1^{j+1} Z_2^k U^l}$ does not appear in the $"\cdots"$ part. Thus, by setting $V_{Z_1^{j+1} Z_2^k U^l}=0$, one plainly solves this recurrence relation to normalize $\alpha_{Z_1^{j+1} Z_2^k U^l}$.
\endproof

\begin{Remark}
\label{rem-orig-pres}
The selected order zero cross-section $Z_1=Z_2=U=V=0$ forces the final normal form transformation of $M^5$ to being origin-preserving.
\end{Remark}

\subsection{Order two}
In this order, we shall consider the recurrence relations
\begin{equation}\label{rec-rel-ord-2}
\aligned
dV_{Z_1\oZ_1}&=\varpi_{Z_1\oZ_1}-V_{Z_1\oZ_1}\big(\mu^1_{Z_1}+\omu^1_{\oZ_1}-\alpha_U\big)-V_{Z_1\oZ_2}\,\omu^2_{\oZ_1}-V_{Z_2\oZ_1} \,\mu^2_{Z_1},
\\
dV_{Z_1\oZ_2}&=\varpi_{Z_1\oZ_2}-V_{Z_1\oZ_1}\,\omu^1_{\oZ_2}-V_{Z_2\oZ_2}\,\mu^2_{Z_1}-V_{Z_1\oZ_2}\big(\mu^1_{Z_1}+\omu^2_{\oZ_2}-\alpha_U\big),
\\
dV_{Z_2\oZ_2}&=\varpi_{Z_2\oZ_2}-V_{Z_1\oZ_2}\,\mu^1_{Z_2}-V_{Z_2\oZ_1}\,\omu^1_{\oZ_2}-V_{Z_2\oZ_2}\big(\mu^2_{Z_2}+\omu^2_{\oZ_2}-\alpha_U\big).
\endaligned
\end{equation}
These relations will be of use in the normalization process whenever at least one of the involved invariants $V_{Z_1\oZ_1}, V_{Z_1\oZ_2}, V_{Z_2\oZ_2}$ was nonzero. But, unfortunately, our assumption that $\bp$ is of Levi non-uniform rank zero under holomorphic transformations prevents us to benefit the application of these recurrence relations.

\begin{Remark}
Denoting
\[
\Delta=V_{Z_1\oZ_1} V_{Z_2\oZ_2}-V_{Z_1\oZ_2}V_{Z_2\oZ_1},
\]
and after applying necessary computations, the above recurrence relations \eqref{rec-rel-ord-2} give that
\[
 d\Delta=\big(2\,\alpha_U-\mu^1_{Z_1}-\mu^2_{Z_2}-\omu^1_{\oZ_1}-\omu^2_{\oZ_2}\big)\,\Delta.
\]
Thus, $\Delta$ is a {\it relative (lifted) invariant} (\cite{Olver-Fels-97}) defined on the lifted bundle $\mathcal B^{(2)}$. Since at each step of normalizations, the value of $\Delta$ at $\bp$ corresponds to the determinant of the Levi matrix \eqref{Levi-matrix-p}, then its relative invariancy confirms the well-known fact that the degeneracy remains invariant under holomorphic transformations (cf. \cite[Proposition 3.6]{Olver-Fels-97}). Let $\mathcal G^{\sf red}$ be the reduction of the holomorphic pseudo-group $\mathcal G$ under the normalizations applied so far. Restricting the recurrence relations \eqref{rec-rel-ord-2} on the fiber $\mathcal G^{\sf red}\times\{\bp\}$, simply gives that $dV_{Z_j\oZ_k}=\varpi_{Z_j\oZ_k}$ for $j,k=1,2$. Thus on this fiber, the three lifted invariants $V_{Z_1\oZ_1}, V_{Z_1\oZ_2}, V_{Z_2\oZ_2}$ are independent of the group parameters. Consequently, the property of $\bp$ having the Levi non-uniform rank zero  remains invariant under the subsequent holomorphic normalizations.
\end{Remark}

\medskip
\noindent
{\bf Notation.} Henceforth and in each step of normalizations, we denote by $\mathcal B_{\bp}$ the fiber $\mathcal G^{\sf red}\times\{\bp\}\subset\mathcal B^{(\infty)}$ over the point $\bp$. As we will see, most of our computations will be done in this domain.

\subsection{Order three}

 While in the former order two, the Levi non-uniform rank zero of $\bp$ played the key role, in this order the $2$-nondegeneracy assumption is essentially effective.
By \cite[eq. (7.1.3--5)]{Ebenfelt-98}, the a partially normalized hypersurface $v=v(z,\oz,u)$ is $2$-nondegenerate whenever its corresponding order three jets enjoy
\begin{equation}\label{2-nondeg-vectors}
{\rm Span}\big\{(v_{z_1^2\oz_1}, v_{z_1^2\oz_2}), (v_{z_1z_2\oz_1}, v_{z_1z_2\oz_2}), (v_{z_2^2\oz_1}, v_{z_2^2\oz_2})\big\}|_{\bp}=\mathbb C^2.
\end{equation}
Corresponding to these jet variables we have, modulo the horizontal coframe, the following six recurrence relations
\begin{equation}\label{rec-rel-ord-3-first-part}
\aligned
dV_{Z_1^2\oZ_1}&\equiv V_{Z_1\oZ_1}\big(4{\rm i}\,V_{Z_1\oZ_2}\omu^2_{U}-\mu^1_{Z_1^2}\big)+4{\rm i}\,V_{Z_1\oZ_1}^2\omu^1_U-V_{Z_1^2\oZ_2}\, \omu^2_{\oZ_1}-2\,V_{Z_1Z_2\oZ_1}\, \mu^2_{Z_1}-V_{Z_2\oZ_1} \mu^2_{Z_1^2}
\\
&+V_{Z_1^2\oZ_1}\big(\alpha_U-\omu^1_{\oZ_1}-2\mu^1_{Z_1}\big),
\\
dV_{Z_1^2\oZ_2}&\equiv V_{Z_1\oZ_2}\big(4{\rm i}\,V_{Z_1\oZ_1}\omu^1_U-\mu^1_{Z_1^2}\big)+4{\rm i}\,V_{Z_1\oZ_2}^2\omu^2_U-V_{Z_1^2\oZ_1}\omu^1_{\oZ_2}-2\,V_{Z_1Z_2\oZ_2}\mu^2_{Z_1}-V_{Z_2\oZ_2}\mu^2_{Z_1^2}
\\
&+V_{Z_1^2\oZ_2} \big(\alpha_U-\omu^2_{\oZ_2}-2\,\mu^1_{Z_1}\big),
\\
dV_{Z_1Z_2\oZ_1}&\equiv V_{Z_1Z_2\oZ_1}\big(\alpha_U-\omu^1_{\oZ_1}-\mu^2_{Z_2}-\mu^1_{Z_1}\big)+V_{Z_1\oZ_1}\big(4{\rm i}\,V_{Z_2\oZ_1}\omu^1_U+2{\rm i}\,V_{Z_2\oZ_2} \omu^2_U-\mu^1_{Z_1Z_2}\big)
\\
&+V_{Z_2\oZ_1} \big(2{\rm i}\,V_{Z_1\oZ_2}\omu^2_U-\mu^2_{Z_1Z_2}\big)-V_{Z_1^2\oZ_1} \mu^1_{Z_2}-V_{Z_1Z_2\oZ_2} \omu^2_{\oZ_1}-V_{Z_2^2\oZ_1} \mu^2_{Z_1},
\\
dV_{Z_1Z_2\oZ_2}&\equiv V_{Z_1Z_2\oZ_2}\big(\alpha_U-\mu^2_{Z_2}-\omu^2_{\oZ_2}-\mu^1_{Z_1}\big)+V_{Z_1\oZ_2} \big(2{\rm i}\,V_{Z_2\oZ_1} \omu^1_U+4{\rm i}\,V_{Z_2\oZ_2} \omu^2_U-\mu^1_{Z_1Z_2}\big)
\\
&+V_{Z_2\oZ_2} \big(2{\rm i}\,V_{Z_1\oZ_1} \omu^1_U-\mu^2_{Z_1Z_2}\big)-V_{Z_1^2\oZ_2} \mu^1_{Z_2}-V_{Z_1Z_2\oZ_1} \omu^1_{\oZ_2}-V_{Z_2^2\oZ_2} \mu^2_{Z_1},
\\
dV_{Z_2^2\oZ_1}&\equiv V_{Z_2^2\oZ_1} \big(\alpha_U-\omu^1_{\oZ_1}-2\,\mu^2_{Z_2}\big)+V_{Z_2\oZ_1} \big(4{\rm i}\,V_{Z_2\oZ_2} \omu^2_U-\mu^2_{Z_2^2}\big)+4{\rm i}\,V_{Z_2\oZ_1}^2\omu^1_U-2\,V_{Z_1Z_2\oZ_1} \mu^1_{Z_2}
\\
&-V_{Z_2^2\oZ_2} \omu^2_{\oZ_1}-V_{Z_1\oZ_1} \mu^1_{Z_2^2}
\\
dV_{Z_2^2\oZ_2}&\equiv V_{Z_2\oZ_2} \big(4{\rm i}\,V_{Z_2\oZ_1} \omu^1_U-\mu^2_{Z_2^2}\big)+V_{Z_2^2\oZ_2} \big(\alpha_U-\omu^2_{\oZ_2}-2\,\mu^2_{Z_2}\big)+4{\rm i}\,V_{Z_2\oZ_2}^2 \omu^2_U -2\,V_{Z_1Z_2\oZ_2} \mu^1_{Z_2}
\\
&-V_{Z_2^2\oZ_1} \omu^1_{\oZ_2} -V_{Z_1\oZ_2} \mu^1_{Z_2^2}.
\endaligned
\end{equation}

\noindent
{\bf Convention.} Hereafter, we use the symbol $"\equiv"$ instead of $"="$, when presenting a recurrence relation modulo the horizontal coframe.
\medskip

In order to normalize the Maurer-Cartan forms by means of the above recurrence relations, we need to know either of the corresponding six differential invariants is nonzero. This necessitate information derived from the $2$-nondegeneracy condition. One verifies that  \eqref{2-nondeg-vectors} is satisfied if at least one of the following combinations is nonzero at $\bp$ (cf. \cite[eq. (7.1.19)]{Ebenfelt-98})
\begin{equation}\label{2-nondeg-cond-2}
\aligned
\Delta_{12}&:=v_{z_1z_2\oz_1}\cdot v_{z_1^2\oz_2}-v_{z_1z_2\oz_2}\cdot v_{z_1^2\oz_1},
\\
\Delta_{23}&:=v_{z_1z_2\oz_2}\cdot v_{z_2^2\oz_1}-v_{z_1z_2\oz_1}\cdot v_{z_2^2\oz_2},
\\
\Delta_{13}&:=v_{z_1^2\oz_1}\cdot v_{z_2^2\oz_2}-v_{z_1^2\oz_2}\cdot v_{z_2^2\oz_1}.
\endaligned
\end{equation}

We emphasize that if two of the above combinations $\Delta_{12}, \Delta_{23}, \Delta_{13}$ vanish at $\bp$ then the third combination has to remain nonzero at this point under holomorphic transformations. We prove this assertion in the case where $\Delta_{23}(\bp)=\Delta_{13}(\bp)=0$ but $\Delta_{12}(\bp)\neq 0$. The proof of the other possible cases is similar.

\begin{Proposition}
\label{propo-invariancy-2-nondeg}
(cf. \cite[Assertion 7.2.1]{Ebenfelt-98}).
Let $\Delta_{23}$ and $\Delta_{13}$ vanish identically at $\bp=0$. Then the value of $\Delta_{12}$ at this point does not vanish under holomorphic transformations.
\end{Proposition}

\proof
By abuse of notation, we continue to write $\Delta_{12}$ for the lift of the above combination $\Delta_{12}$ under holomorphic transformations, i.e.
\[
\Delta_{12}:=V_{Z_1Z_2\oZ_1}\cdot V_{Z_1^2\oZ_2}-V_{Z_1Z_2\oZ_2}\cdot V_{Z_1^2\oZ_1}.
\]
Taking into account that the order two lifted invariants $V_J$, $\#J=2$ vanish identically on the fiber $\mathcal B_{\bp}$ over $\bp$, the recurrence relations \eqref{rec-rel-ord-3-first-part} give
\[
\aligned
d\Delta_{12}\equiv\big(2\,\alpha_U-3\,\mu^1_{Z_1}-\mu^2_{Z_2}-\omu^1_{\oZ_1}-\omu^2_{\oZ_2}\big) \Delta_{12}.
\endaligned
\]
Thus on the mentioned fibre, the pseudo-group acts by scaling on $\Delta_{12}$. Hence, it remains nonzero on this bundle under holomorphic transformations.
\endproof

The $2$-nondegeneracy condition \eqref{2-nondeg-cond-2} provides various possibilities of normalizations in order three, each of which produces a certain branch in the normalization process. Ebenfelt realized in \cite{Ebenfelt-98} these branches as \eqref{Ebenfelt-Branches}. Our main goal in this paper is to complete the normalizations in each branch and construct their associated complete normal forms.

\begin{Remark}
\label{rem-r-neq-1-2}
Ebenfelt claims in \cite[p. 339]{Ebenfelt-98} that the lifted differential invariant $V_{Z_1Z_2\oZ_2}$ which corresponds to $b_2$ in \cite{Ebenfelt-98}, can be normalized to zero in all branches. As we will see, this claim is true except in a very specific case. Indeed, when we have on the fiber $\mathcal B_{\bp}$ that $|V_{Z_1^2\oZ_2}|=|V_{Z_1Z_2\oZ_1}|$ and $V_{Z_2^2\oZ_2}=0$\,\,---\,\,this can occur in branch (A.ii.1) of \eqref{Ebenfelt-Branches}\,\,---\,\,then setting $V_{Z_1Z_2\oZ_2}=0$ in the fourth equation in \eqref{rec-rel-ord-3-first-part} is superfluous as in this case only real or imaginary part of the Maurer-Cartan form $\mu^1_{Z_2}$ is normalizable. This situation corresponds to the specific case $|r|=\frac{1}{2}$ where $r$ is the invariant defined in \cite[eq. (7.2.13)]{Ebenfelt-98} (this exception was also realized by Ershova in \cite[p. 191]{Ershova-01}). Nevertheless, in this paper we do not aim to consider this specific case for two reasons. First, we do not like to make the paper longer and second, one can get rid of this phenomenon by appropriate (but different from what we do in the next section) normalizations of the two lifted invariants $V_{Z_1^2\oZ_2}$ and $V_{Z_1Z_2\oZ_1}$ in branch (A.ii.1).
\end{Remark}

Let us conclude this section by noticing that in the present order three, we have in addition the three recurrence relations
\begin{equation}\label{rec-rel-ord-3-compl}
\aligned
dV_{Z_1\oZ_1U}&\equiv V_{Z_1\oZ_1} \big(\alpha_{UU}-\mu^1_{Z_1U}-\omu^1_{\oZ_1U}\big)-V_{Z_1\oZ_1 U} \big(\mu^1_{Z_1}+\omu^1_{\oZ_1}\big)-V_{Z_1^2\oZ_1} \mu^1_U-V_{Z_1\oZ_1^2} \omu^1_U
\\
&-
V_{Z_2\oZ_1} \mu^2_{Z_1U}-V_{Z_1\oZ_2} \omu^2_{\oZ_1U}-V_{Z_2\oZ_1U} \mu^2_{Z_1}-V_{Z_1\oZ_2U} \omu^2_{\oZ_1}-V_{Z_1Z_2\oZ_1} \mu^2_U-V_{Z_1\oZ_1\oZ_2} \omu^2_U,
\\
dV_{Z_1\oZ_2U}&\equiv V_{Z_1\oZ_2} \big(\alpha_{UU}-\mu^1_{Z_1U}-\omu^2_{\oZ_2U}\big)-V_{Z_1\oZ_2U} \big(\mu^1_{Z_1}+\omu^2_{\oZ_2}\big)-V_{Z_1\oZ_1} \omu^1_{\oZ_2U}-V_{Z_2\oZ_2} \mu^2_{Z_1U}
\\
&-V_{Z_1^2\oZ_2} \mu^1_U-V_{Z_1\oZ_1U}\omu^1_{\oZ_2}-V_{Z_2\oZ_2U} \mu^2_{Z_1}-V_{Z_1Z_2\oZ_2} \mu^2_U-V_{Z_1\oZ_1\oZ_2}\omu^1_U-V_{Z_1\oZ_2^2} \omu^2_U,
\\
dV_{Z_2\oZ_2U}&\equiv V_{Z_2\oZ_2} \big(\alpha_{UU}-\mu^2_{Z_2U}-\omu^2_{\oZ_2U}\big)-V_{Z_2\oZ_2U} \big(\mu^2_{Z_2}+\omu^2_{\oZ_2}\big)-V_{Z_1\oZ_2} \mu^1_{Z_2U}-V_{Z_2\oZ_1} \omu^1_{\oZ_2U}
\\
&- V_{Z_1Z_2\oZ_2} \mu^1_{U}-V_{Z_2\oZ_1\oZ_2} \omu^1_U-V_{Z_1\oZ_2U} \mu^1_{Z_2} -V_{Z_2\oZ_1U} \omu^1_{\oZ_2} -V_{Z_2^2\oZ_2} \mu^2_U - V_{Z_2\oZ_2^2} \omu^2_U.
\endaligned
\end{equation}

In the following three sections, our objective is to complete the partial normal forms \eqref{Ebenfelt-Branches} across the five emerged branches.

\section{Branches (A.{\rm ii}.1) and (A.{\rm ii}.2)}
\label{sec-A.ii.1-A.ii.2}

Since both the first two partial normal forms (A.ii.1) and (A.ii.2) in \eqref{Ebenfelt-Branches} fall under the situation where the value of $\Delta_{12}$ is nonzero at the point $\bp = 0$, we study them together. The coefficient of the monomial $z_1z_2\oz_1$ is nonzero in these branches and thus, we are permitted here to set the corresponding lifted invariant $V_{Z_1Z_2\oZ_1}$ to some nonzero constant number, say
\[
V_{Z_1Z_2\oZ_1}=1.
\]
Accordingly, one can solve the third recurrence relation in \eqref{rec-rel-ord-3-first-part} to normalize
\[
\aligned
\mu^2_{Z_2}&\equiv V_{Z_1\oZ_1} \big(4{\rm i}\,V_{Z_2\oZ_1} \omu^1_U+2{\rm i}\,V_{Z_2\oZ_2} \omu^2_U- \mu^1_{Z_1Z_2}\big)+V_{Z_2\oZ_1} \big(2{\rm i}\,V_{Z_1\oZ_2} \omu^2_U- \mu^2_{Z_1Z_2}\big)
\\
&-V_{Z_1^2\oZ_1} \mu^1_{\oZ_2}-V_{Z_2^2\oZ_1} \mu^2_{Z_1}+\alpha_U-\mu^1_{Z_1}-\omu^1_{\oZ_1}.
\endaligned
\]
Even more generally, we have

\begin{Lemma}
\label{mu2-Z1-Z2-Zbar1}
For every $j,k,l\geq 0$ with $(j,l,k)\neq (0,0,0)$, solving the recurrence relation of $V_{Z_1^{j+1} Z_2^{k+1}\oZ_1U^l}=0$ provides the normalization of the Maurer-Cartan form $\mu^2_{Z_1^jZ_2^{k+1}U^l}$.
\end{Lemma}

\proof
By analysing the vector component $\phi^{z_1z_2\oz_1}$ of the prolonged vector field \eqref{v-infty}, one finds that
\[
\phi^{z_1z_2\oz_1}=-v_{z_1z_2\oz_1}\,\xi^2_{z_2}+\cdots
\]
where $"\cdots"$ stands for terms that do not involve the vector component $\xi^2_{z_2}$. Then, the prolongation formula \eqref{eq: prolongation formula 2} results for each $j, k ,l\geq 0$ that
\[
\phi^{z_1^{j+1}z_2^{k+1}\oz_1u^l}=-v_{z_1z_2\oz_1}\,\xi^2_{z_1^jz_2^{k+1}u^l}+\cdots
\]
where $\xi^2_{z_1^jz_2^{k+1}u^l}$ does not appear in the part $"\cdots"$. Then, in light of the recurrence formula \eqref{rec-formula} and after invariantization, we have
\[
dV_{Z_1^{j+1}Z_2^{k+1}\oZ_1U^l}\equiv -V_{Z_1Z_2\oZ_1}\,\mu^2_{Z_1^jZ_2^{k+1}U^l}+\cdots.
\]
Thus, with $V_{Z_1Z_2\oZ_1}=1$ and $V_{Z_1^{j+1}Z_2^{k+1}\oZ_1U^l}=0$, one may solve the above recurrence relation to normalize the Maurer-Cartan form $\mu^2_{Z_1^jZ_2^{k+1}U^l}$.
\endproof

Furthermore, with the assumption\footnote{As mentioned in Remark \ref{rem-r-neq-1-2}, we do not aim to consider the specific opposite case.} $|V_{Z_1^2\oZ_2}|\neq|V_{Z_1Z_2\oZ_1}|$, the recurrence relation of $V_{Z_1Z_2\oZ_2}$ enables one to normalize the entire complex Maurer-Cartan form $\mu^1_{Z_2}$ if we set\footnote{Henceforth and due to their length, we do not present the normalized expressions of the Maurer-Cartan forms. However, they are available in the {\sc Maple} worksheet \cite{Maple}.}
\[
V_{Z_1Z_2\oZ_2}=0.
\]
Proceeding along the same arguments of the proofs of Lemmas \ref{lem-ord-one} and \ref{mu2-Z1-Z2-Zbar1} and by the careful analysis of the prolongation formula \eqref{eq: determining equations}, one finds in general that

\begin{Lemma}
\label{mu1-Z1-Z2-Zbar2}
For every $j,k,l\geq 0$, one can normalize the Maurer-Cartan form $\mu^1_{Z_1^jZ_2^{k+1}U^l}$ by specifying $V_{Z_1^{j+1}Z_2^{k+1}\oZ_2U^l}=0$.
\end{Lemma}

For the next step of normlizations, we consider the lifted invariant $V_{Z_1^2\oZ_2}$ in which its recurrence relation in \eqref{rec-rel-ord-3-first-part} is now of the form
\[
dV_{Z_1^2\oZ_2}\equiv -2{\rm i}\,V_{Z_1^2\oZ_2}\, {\rm Im}\mu^1_{Z_1}+\cdots
\]
where the $"\cdots"$ part stands for the terms that vanish identically on the fiber $\mathcal B_{\bp}$. Setting ${\rm Im}V_{Z_1^2\oZ_2}=0$ enables one to solve the imaginary part of the above relation to normalize the real Maurer-Cartan form ${\rm Im}\mu^1_{Z_1}$. More generally, we have

\begin{Lemma}
\label{lem-Im-mu1-Z1-U}
Let $j, l\geq 0$. The real maurer-Cartan form ${\rm Im}\mu^1_{Z_1^{j+1}U^l}$ can be normalized by setting ${\rm Im}V_{Z_1^{j+2}\oZ_2U^l}=0$.
\end{Lemma}

On the fiber $\mathcal B_{\bp}$ and modulo the horizontal coframe, the recurrence relation of the real differential invariant ${\rm Re}V_{Z_1^2\oZ_2}$ is in turn
 \[
 d{\rm Re}V_{Z_1^2\oZ_2}\equiv 0.
 \]
Accordingly, when restricted to $\mathcal B_{\bp}$, we have ${\rm Re}V_{Z_1^2\oZ_2}$ independent of the holomorphic group parameters. Thus, further normalizations will not effect it in branches (A.ii.1) and (A.ii.2).
Let us denote
\[
r:=\frac{1}{2}\,{\rm Re}V_{Z_1^2\oZ_2}(\bp).
\]
Here, the coefficient $\frac{1}{2}$ provides the consistency of this notation with Ebenfelt's invariant introduced in \cite[eq. (7.2.13)]{Ebenfelt-98}. By Proposition \ref{propo-invariancy-2-nondeg} and  after setting ${\rm Im}V_{Z_1^2\oZ_2}=0$, we have $r\neq 0$.

\begin{Remark}
We can regard $r$ as a positive integer. Indeed, even if $r<0$, the simple transformation $z_1\mapsto {\rm i} z_1$ converts the coefficient $r$ of the monomial $z_1^2\oz_2$ in the defining function of our partial normal form to $-r$, while of the coefficient of $z_1z_2\oz_1$ remains $1$, unaffected by the transformation.
\end{Remark}

The next candidate for normalization is $V_{Z_1^2\oZ_1}$. After applying the previous normalizations and on the fiber $\mathcal B_{\bp}$, the recurrence relation of this invariant takes the form
\begin{equation}\label{dV-z1z1oz1}
\aligned
dV_{Z_1^2\oZ_1}&\equiv V_{Z_1^2\oZ_1}\,\big(\alpha_U-3\,{\rm Re}\mu^1_{Z_1}+\cdots\big)+\frac{1}{8r}\,\big(2r\,V_{Z_1^2\oZ_1}V_{Z_2^2\oZ_1} - V_{Z_1^2\oZ_1}^2V_{Z_2^2\oZ_2} - 16\,r\big)\,\mu^2_{Z_1}
\\
+&\frac{1}{8r}\,\big(-2r\,V_{Z_1\oZ_2^2}V_{Z_1^2\oZ_1} + V_{Z_1\oZ_1^2}V_{Z_2\oZ_2^2}V_{Z_1^2\oZ_1}-16\,r^2\big)\,\omu^2_{\oZ_1}+\cdots.
\endaligned
\end{equation}
The feasibility of practical normalization by means of this relation depends on whether we have $r\neq 1$ or $r=1$. In the former case, one can readily normalize $\mu^2_{Z_1}$ by setting $V_{Z_1^2\oZ_1}=0$. But, when $r=1$ this approach is not applicable. In this case, the application of the above relation relies upon vanishing/nonvanishing of $V_{Z_1^2\oZ_1}$. Accordingly, we shall split the next computations into three branches
\begin{itemize}
  \item[{\rm (A$'$.ii.1)}.] when $r\neq 1$.
  \item[{\rm (A$''$.ii.1)}.] when $r=1$ and $V_{Z_1^2\oZ_1}=0$.
  \item[{\rm (A.ii.2)}.] when $r=1$ and $V_{Z_1^2\oZ_1}\neq0$.
\end{itemize}
Actually, (A$'$.ii.1) and (A$''$.ii.1) divide the branch (A.ii.1) into two parts. Moreover, (A.ii.2) coincides with its correspondence in \eqref{Ebenfelt-Branches}.
 Before proceeding into the next computations along the above three branches, we notice that in light of the normalizations performed thus far, the remaining unnormalized
Maurer–Cartan forms are
 \begin{equation}\label{MC-2}
 \mu^1_{U^{l+1}}, \qquad \omu^1_{U^{l+1}}, \qquad {\rm Re}\mu^1_{Z_1^{j+1}U^l}, \qquad \mu^2_{Z_1^jU^l}, \qquad \omu^2_{\oZ_1^jU^l}, \qquad \alpha_{U^{l+1}}
 \end{equation}
 for $j, l\geq 0$.

\subsection{Branch (A$'$.ii.1)}

Thanks to the assumption $r\neq 1$, plainly setting
\[
V_{Z_1^2\oZ_1}=0
\]
provides the opportunity of normalizing the complex Maurer-Cartan form $\mu^2_{Z_1}$ on the fiber $\mathcal B_{\bp}$\,\,---\,\,and thus in a neighborhood of it\,\,---\,\,by solving the greatly simplified recurrence relation (cf.  \eqref{dV-z1z1oz1})
\begin{equation*}
0=dV_{Z_1^2\oZ_1}\equiv 4{\rm i}\,V_{Z_1\oZ_1}^2\,\omu^1_U + V_{Z_1\oZ_1}\,\big(-4{\rm i}\,V_{Z_1\oZ_2}\,\omu^2_U - \omu^1_{Z_1^2}\big) - V_{Z_2\oZ_1}\,\mu^2_{Z_1^2}- 2r\,\omu^2_{\oZ_1}- 2\,\mu^2_{Z_1}.
\end{equation*}
More generally we have

\begin{Lemma}
\label{lem-mu2-Z1}
In branch {\rm (A$'$.ii.1)} and for each $j,l\geq 0$, one normalizes the Maurer-Cartan form $\mu^2_{Z_1^{j+1}U^l}$ by setting $V_{Z_1^{j+2}\oZ_1U^l}=0$.
\end{Lemma}

Now let us inspect the last two recurrence relations in \eqref{rec-rel-ord-3-first-part} which are not considered yet. Our computations \cite{Maple} show that on the fiber $\mathcal B_{\bp}$ we have them simply as
\[
\aligned
dV_{Z_2^2\oZ_1}&\equiv \big(3\,{\rm Re}\mu^1_{Z_1}-\alpha_U\big)\,V_{Z_2^2\oZ_1},
\\
dV_{Z_2^2\oZ_2}&\equiv 2\,\big(3\,{\rm Re}\mu^1_{Z_1}-\alpha_U\big)\,V_{Z_2^2\oZ_1}.
\endaligned
\]
Thus, the partially normalized pseudo-group acts by scaling on $V_{Z_2^2\oZ_1}$ and $V_{Z_2^2\oZ_1}$, when restricted to the fiber $\mathcal B_{\bp}$. Since, the coefficients of the corresponding monomials $z_2^2\oz_1$ and $z_2^2\oz_2$ are zero in the initial expression (A.ii.1) in \eqref{Ebenfelt-Branches}, then they remain vanished under holomorphic transformations. Consequently, here on the fiber $\mathcal B_{\bp}$ we have
\[
V_{Z_2^2\oZ_1}=V_{Z_2^2\oZ_2}=0.
\]

Let us consider at the three remained order three recurrence relations \eqref{rec-rel-ord-3-compl}. Our computations show that at this stage of normalizations, the first relation is now of the form
\[
\aligned
dV_{Z_1\oZ_1U}&\equiv\big(\frac{2{\rm i}\,(r\,V_{Z_2\oZ_1U} - V_{Z_1\oZ_2U}) V_{Z_2\oZ_1}V_{Z_1\oZ_1}}{r^2-1}-1\big)\,\mu^2_U
\\
&+\big(\frac{-2{\rm i}\, (rV_{Z_1\oZ_2U} - V_{Z_2\oZ_1U}) V_{Z_1\oZ_2}V_{Z_1\oZ_1}}{r^2-1}-1\big)\,\omu^2_U+\cdots.
\endaligned
\]
On the fiber $\mathcal B_{\bp}$, the written part at the right hand side of the above relation simplifies to $-2{\rm Re}\mu^2_U$. Thus, provided
\[
V_{Z_1\oZ_1U}=0,
\]
 it is possible in a local neighborhood of this fiber to normalize the real Maurer-Cartan form ${\rm Re}\mu^2_U$, by solving the above recurrence relation.
In general we have

\begin{Lemma}
\label{lem-Re-mu2-U}
For $l\geq 0$ and in branch {\rm (A$'$.ii.1)}, solving the recurrence relation of $V_{Z_1\oZ_1U^{l+1}}=0$ offers the normalization of the real Maurer-Cartan form ${\rm Re}\mu^2_{U^{l+1}}$.
\end{Lemma}

Next, let us consider the second recurrence relation $dV_{Z_1\oZ_2U}$ in \eqref{rec-rel-ord-3-compl}. On the fiber $\mathcal B_{\bp}$, its too large expression can be written as
\[
dV_{Z_1\oZ_2U}\equiv -2r\,\mu^1_U-\omu^1_U+\cdots
\]
where $"\cdots"$ involves no nonzero coefficient of $\mu^1_U$ or its conjugation. Thus, reminding that we excluded the specific case $r=\frac{1}{2}$, this recurrence relation offers to normalize the Maurer-Cartan form $\mu^1_U$ by setting
\[
V_{Z_1\oZ_2U}=0.
\]
More generally, we have

\begin{Lemma}
\label{lem-mu1-U}
In branch {\rm (A$'$.ii.1)} and for each $l\geq 0$, one normalizes the Maurer-Cartan form $\mu^1_{U^{l+1}}$ by setting $V_{Z_1\oZ_2U^{l+1}}=0$.
\end{Lemma}

At this stage, the list of the remained yet unnormalized Maurer-Cartan forms \eqref{MC-2} is reduced to the real forms
 \begin{equation}\label{MC-3-A-1}
{\rm Re}\mu^1_{Z_1^{j+1}U^l}, \qquad {\rm Im}\mu^2_{U^{l+1}}, \qquad \alpha_{U^{l+1}}
 \end{equation}
for $j,l\geq 0$. In the current order three, none of the lifted differential invariant has the potential of normalizing these forms. Then, we have to proceed the computations into the next order four.

\subsubsection{Order four of branch {\rm (A$'$.ii.1)}}

Our tedious computations of all recurrence relations in order four revealed that only one lifted differential invariant, namely ${\rm Re}V_{Z_1^3\oZ_2}$, is of help in the process of normalization. Restricted to the fiber $\mathcal B_{\bp}$, the recurrence relation of this invariant is
\begin{equation}\label{dV-Z1-3-Zbar2}
d{\rm Re}V_{Z_1^3\oZ_2}\equiv -6r\,{\rm Re}\mu^1_{Z_1^2}+\cdots.
\end{equation}
This equation suggests to normalize the real Maurer-Cartan form ${\rm Re}\mu^1_{Z_1^2}$ by solving it after setting
\[
{\rm Re}V_{Z_1^3\oZ_2}=0.
\]
More generally, we have

\begin{Lemma}
\label{lem-Re-mu1-Z1Z1}
For every $j, l\geq 0$ and in branch {\rm (A$'$.ii.1)}, one can normalize the real Maurer-Cartan form ${\rm Re}\mu^1_{Z_1^{j+2}U^l}$ by specifying ${\rm Re}V_{Z_1^{j+3}\oZ_2U^l}=0$.
\end{Lemma}

Combining the above observation with Lemma \ref{lem-Im-mu1-Z1-U}, it becomes evident setting $V_{Z_1^{3+j}\oZ_2U^l}=0$ is sufficient for normalizing the entire complex Maurer-Cartan form $\mu^1_{Z_1^{j+2}U^l}$.

Some of the other order four recurrence relations are applicable in the normalization process only when they does not vanish. For the sake of generality and to prevent ourselves from producing further subbranches, let us neglect the contribution of such recurrence relations and proceed into the next order with the hope of finding further general normalizations.
Notice that at this stage, the list of the yet unnormalized Maurer-Cartan forms \eqref{MC-3-A-1} is reduced to
  \begin{equation}\label{MC-4-A-1}
{\rm Re}\mu^1_{Z_1U^l}, \qquad {\rm Im}\mu^2_{U^{l+1}}, \qquad \alpha_{U^{l+1}}, \qquad  l\geq 0.
 \end{equation}

\subsubsection{Order five of branch {\rm (A$'$.ii.1)}}

We start this order by the lifted invariant $V_{Z_1^2\oZ_1^2\oZ_2}$ which, on the fiber $\mathcal B_{\bp}$, has the recurrence relation
\begin{equation}\label{dV-Z1-2-Zbar1-2-Zbar2}
dV_{Z_1^2\oZ_1^2\oZ_2}\equiv 4{\rm i}\,\big(\omu^2_U-(1+r^2)\mu^2_U)+\cdots.
\end{equation}
This relation provides the normalization of ${\rm Im}\mu^2_U$ by setting ${\rm Re}V_{Z_1^2\oZ_1^2\oZ_2}=0$. In general, we have

\begin{Lemma}
\label{lem-Im-mu2-U}
For every $l\geq 0$ and in branch {\rm (A$'$.ii.1)} the recurrence relation of ${\rm Re}V_{Z_1^2\oZ_1^2\oZ_2U^l}=0$ can be solved to normalize the real Maurer-Cartan form ${\rm Im}\mu^2_{U^{l+1}}$.
\end{Lemma}

Similar to the former order four, further normalizations of the Maurer-Cartan forms are available in this order only if certain fifth order differential invariants does not vanish. Again for the desire of mere generality, we choose not involve ourselves with these possibilities and move on to the next order six to seek further normalizations.
Notice that the collection of the yet unnormalized Maurer-Cartan forms \eqref{MC-4-A-1} is now reduced to

  \begin{equation}\label{MC-5-A-1}
{\rm Re}\mu^1_{Z_1U^l}, \qquad \alpha_{U^{l+1}}, \qquad {\rm for} \ l\geq 0.
 \end{equation}

\subsubsection{Order six of branch {\rm (A$'$.ii.1)}}

Contrary to the previous two orders, it appears in order six a plenty of differential invariants with the potential of normalizing Maurer-Cartan forms. Roughly speaking, the main reason of this change in the behaviour of differential invariants is the assumption, made regarding the rank of the Levi matrix of $M^5$ at point $\bp$, which was assumed to be zero. It necessitates the lower order invariants to undergo extra differentiations for reaching to order three, where certain lifted invariants are normalized to non-constant integers.

Among the already mentioned differential invariants, we choose to consider the imaginary parts of two of them, namely $V_{Z_1^4Z_2^2}$ and $V_{Z_1^3\oZ_1\oZ_2^2}$, where their recurrence relations on the fiber $\mathcal B_{\bp}$ are
\begin{equation}\label{dImV-Z14-Zbar22}
\aligned
d{\rm Im}V_{Z_1^4\oZ_2^2}&\equiv \big(48\,\alpha_{UU}-192\,{\rm Re}\mu^1_{Z_1U}\big)\,r^2+\cdots,
\\
d{\rm Im}V_{Z_1^3\oZ_1\oZ_2^2}&\equiv \big(24\,\alpha_{UU}-72\,{\rm Re}\mu^1_{Z_1U}\big)\,r+\cdots.
\endaligned
\end{equation}
By setting ${\rm Im}V_{Z_1^4\oZ_2^2}={\rm Im}V_{Z_1^3\oZ_1\oZ_2^2}=0$, one readily solves the above two recurrence relations to normalize the real Maurer-Cartan forms $\alpha_{UU}$ and ${\rm Re}\mu^1_{Z_1U}$. More generally we have

\begin{Lemma}
\label{lem-alpha-UU}
Let $l\geq 0$. Then, in branch {\rm (A$'$.i.1)}
\begin{itemize}
  \item[$1)$] one can normalize the Maurer-Cartan form $\alpha_{U^{2+l}}$ by setting ${\rm Im}V_{Z_1^4\oZ_2^2U^l}=0$.
  \item[$2)$] one can normalize the Maurer-Cartan form ${\rm Re}\mu^1_{Z_1U^{l+1}}$ by setting ${\rm Im}V_{Z_1^3\oZ_1\oZ_2^2U^l}=0$.
\end{itemize}
\end{Lemma}

We have succeeded to normalize all but only two of the remained Maurer–Cartan forms
\[
{\rm Re}\mu^1_Z, \qquad {\rm and} \qquad \alpha_U.
\]
This, according to \cite{Valiquette-SIGMA}, implies that the isotropy group at the point $\bp$ of $2$-nondegenerate hypersurfaces in branch (A$'$.ii.1) are of dimensions $\leq 2$. Thus, at this branch, we have reached Ershova's bound \cite{Ershova-01} for the corresponding isotropy groups. The maximum possible dimension two is enjoyed by those hypersurfaces in which the remained unconsidered lifted invariants $V_J$ vanish identically on them. They actually are Ershova's {\it model hypersurfaces}
\[
v=z_1z_2\oz_1+z_1\oz_1\oz_2+r\,\big(z_1^2\oz_2+z_2\oz_1^2\big), \qquad r\neq 1
\]
which admit the real parts of the {\it dilation} vector fields
\begin{equation}\label{dilation-fields}
{\sf D}_1:=z_1\partial_{z_1}-2\,z_2\partial_{z_2}, \qquad {\sf D_2}:=w\partial_{w}+z_2\partial_{z_2}
\end{equation}
 as the generators of their isotropy algebras.

We are now ready to present the complete normal form of branch (A$'$.ii.1). Recall that for every pair $\ell=(l_1, l_2)\in\mathbb N_0^2$, we denoted $z^\ell=z_1^{l_1}z_2^{l_2}$. Moreover, we let $|\ell|:=l_1+l_2$ and $\ell!:=l_1!\,l_2!$.

 \begin{Theorem}
 Let $M^5\subset\mathbb C^3$ be a $2$-nondegenerate real hypersurface of Levi non-uniform rank zero at the origin point $\bp=0$. If $M^5$ belongs to the branch {\rm (A$'$.ii.1)}, then it can be mapped through an origin-preserving transformation to the complete normal form
 \begin{equation}\label{NF-A'.ii.1}
 \aligned
 v= z_1z_2\oz_1+z_1\oz_1\oz_2+r\,\big(z_1^2\oz_2+z_2\oz_1^2\big)+V_{Z_2\oZ_2U}\,z_2\oz_2u+\sum_{|\ell_1|+|\ell_2|+l\geq 4}\,\frac{V_{Z^{\ell_1}\oZ^{\ell_2}U^l}}{\ell_1!\,\ell_2!\,l!}\,z^{\ell_1}\oz^{\ell_2} u^l,
 \endaligned
 \end{equation}
 for a unique real number $r\neq 1$. Moreover, regarding the conjugation relation, the coefficients $V_J$ enjoy the cross-section
 \[
 \aligned
 0\equiv V_{Z^\ell U^l}&=V_{Z_1^{j+1}Z_2^{k+1}\oZ_tU^l}=V_{Z_1^{j+3}\oZ_2U^l}=V_{Z_1^{j+2}\oZ_1U^l}=V_{Z_1\oZ_tU^{l+1}}
 \\
 &={\rm Re}V_{Z_1^2\oZ_1^2\oZ_2U^l}={\rm Im}V_{Z_1^2\oZ_2U^l}={\rm Im}V_{Z_1^4\oZ_2^2U^l}={\rm Im}V_{Z_1^3\oZ_1\oZ_2^2U^l}
 \endaligned
 \]
 for $t=1,2$, $\ell\in\mathbb N_0$ and $j,k,l\geq 0$. Furthermore, the isotropy group of $M^5$ at $\bp$ is of dimension at most two.
 \end{Theorem}

 \begin{Remark}
 We emphasize that in branch (A$'$.ii.1), there exists an abundance of hypersurfaces with isotropy groups whose dimensions are {\it absolutely} less than two. It depends to vanishing/nonvanishing of some specific differential invariants on the fiber $\mathcal B_{\bp}$. For example, our computations show that if the third order lifted invariant $V_{Z_2\oZ_2U}$ does not vanish\footnote{By evaluating the third recurrence relation in \eqref{rec-rel-ord-3-compl} on $\mathcal B_{\bp}$, one finds that vanishing/nonvanishing of  $V_{Z_2\oZ_2U}$ is invariant on this fiber.} at $\bp$, then one can normalize either of the remained two Maurer-Cartan forms $\alpha_U, {\rm Re}\mu^1_{Z_1}$ by setting $V_{Z_2\oZ_2U}=1$. In that case, the isotropy group of the appearing normal form
  \[
v=z_1z_2\oz_1+z_1\oz_1\oz_2+r\,\big(z_1^2\oz_2+z_2\oz_1^2\big)+z_2\oz_2u+\sum_{|\ell_1|+|\ell_2|+l\geq 4}\,\frac{V_{Z^{\ell_1}\oZ^{\ell_2}U^l}}{\ell_1!\,\ell_2!\,l!}\,z^{\ell_1}\oz^{\ell_2} u^l, \qquad r\neq 1
\]
is of dimension $\leq 1$. In particular, when all appearing lifted invariants $V_J$, $\#J\geq 4$ vanish identically at $\bp$, then the isotropy group is exactly $1$-dimensional generated infinitesimally by
\[
{\sf D}_1+2\,{\sf D}_2=z_1\partial_{z_1}+2w\,\partial_w.
\]
 \end{Remark}

\subsection{Branch (A$''$.ii.1)}

Now, let us consider the second part of the branch (A.ii.1), where in addition to the original assumptions of this branch, we also suppose that after the partial normalizations made up to Lemma \ref{lem-Im-mu1-Z1-U}, we have
\[
 r=1 \qquad  {\rm and} \qquad V_{Z_1^2\oZ_1}=0.
\]

These assumptions will not  disrupt the normalizations achieved in Lemmas \ref{mu2-Z1-Z2-Zbar1}\,--\,\ref{lem-Im-mu1-Z1-U}. Thus, at this stage, we may assume the collection \eqref{MC-2} as the remaining unnormalized Maurer-Cartan forms.

Since inserting $r=1$ in the first recurrence relation $dV_{Z_1^2\oZ_1}$ of the list \eqref{rec-rel-ord-3-first-part} converts the term $- 2r\omu^2_{\oZ_1}- 2\mu^2_{Z_1}$ to $-4{\rm Re}\mu^2_{Z_1}$ then, unfortunately, this relation can not normalize anymore the entire complex Maurer-Cartan form $\mu^2_{Z_1}$ . Indeed, here, it provides us with normalizing only the real part of $\mu^2_{Z_1}$ after setting
\[
{\rm Re}V_{Z_1^2\oZ_1}=0.
\]

Then\,\,---\,\,in contrast to Lemma \ref{lem-mu2-Z1}\,\,---\,\,we have in general

\begin{Lemma}
\label{lem-mu2-Z1-A''.ii.1}
In branch {\rm (A$''$.ii.1)} and for each $j,l\geq 0$, the recurrence relation of ${\rm Re}V_{Z_1^2\oZ_1U^l}=0$ can be solved to normalizes the Maurer-Cartan form ${\rm Re}\mu^2_{Z_1U^l}$.
\end{Lemma}

As in the branch (A$'$.ii.1), the two lifted invariants $V_{Z_2^2\oZ_1}$ and $V_{Z_2^2\oZ_2}$ vanish identically at the point $\bp$ (see the paragraph after Lemma \ref{lem-mu2-Z1}). Similarly, inspecting the recurrence relation of ${\rm Im}V_{Z_1^2\oZ_1}$ shows that our partially normalized pseudo-group acts by scaling on this real lifted invariant, when restricted to the fiber $\mathcal B_{\bp}$. Thus, the value of ${\rm Im}V_{Z_1^2\oZ_1}$ (and hence $V_{Z_1^2\oZ_1}$) remains zero on this bundle under holomorphic transformations.

Again similar to the branch (A$'$.ii.1), here the lifted invariant $V_{Z_1\oZ_1U}$ provides the opportunity of normalizing the real Maurer-Cartan form ${\rm Re}\mu^2_U$. Thus, Lemma \ref{lem-Re-mu2-U} works as well in this branch.
Moreover, when restricted to the fiber $\mathcal B_{\bp}$, the recurrence relation of $V_{Z_1\oZ_2U}$ exhibits the term $-2\mu^1_U-\omu^1_1$ and thus\,\,---\,\,again similar to (A$'$.ii.1)\,\,---\,\,one can normalize the Maurer-Cartan form $\mu^1_U$ by setting $V_{Z_1\oZ_2U}=0$. Consequently, the general Lemma \ref{lem-mu1-U} holds also in this case.

The only remained third order lifted invariant $V_{Z_2\oZ_2U}$ is of no use to normalize further Maurer-Cartan forms unless it does not vanish. With the aim of seeking for the generality, we do not aim to consider this possibility.

Summing up, then the list of the yet unnormalized Maurer-Cartan forms \eqref{MC-2} is now reduced to
 \begin{equation}\label{MC-3-A''}
 {\rm Re}\mu^1_{Z_1^{j+1}U^l}, \qquad \mu^2_{Z_1^{j+2}U^l}, \qquad {\rm Im}\mu^2_{Z_1U^l}, \qquad {\rm Im}\mu^2_{U^{l+1}} \qquad \alpha_{U^{l+1}}
 \end{equation}
 for $j, l\geq 0$. Let us continue the normalizations in the next order four.

 \subsubsection{Order four of branch {\rm (A$''$.ii.1)}}

Consider the lifted invariant $V_{Z_1^3\oZ_1}$, with the following recurrence relation on $\mathcal B_{\bp}$
\[
dV_{Z_1^3\oZ_1}\equiv \big(\alpha_U-4\,{\rm Re}\mu^1_{Z_1}\big)V_{Z_1^3\oZ_1}-3\,\mu^2_{Z_1^2}.
\]
Thus, setting $V_{Z_1^3\oZ_1}=0$, it offers the normalization of the Maurer-Cartan form $\mu^2_{Z_1^2}$. More generally we have

\begin{Lemma}
In branch {\rm (A$''$.ii.1)} and for each $j,l\geq 0$, one can normalize the Maurer-Cartan form $\mu^2_{Z_1^{j+2}U^l}$ by setting $V_{Z_1^{j+3}\oZ_1U^l}=0$.
\end{Lemma}

Furthermore, on the fiber $\mathcal B_{\bp}$, the recurrence relation of the lifted invariant $V_{Z_1^3\oZ_2}$ is exactly as \eqref{dV-Z1-3-Zbar2} with $r=1$. Then, Lemma \ref{lem-Re-mu1-Z1Z1} works well also in this branch.
Unfortunately, the other fourth order recurrence relations are of no general effect in the normalization process and we shall proceed into the next order for further possible normalizations.

 \subsubsection{Orders five and six of branch {\rm (A$''$.ii.1)}}

 As in the fifth order of branch (A$'$.ii.1), here only the recurrence relation of $V_{Z_1^2\oZ_1^2\oZ_2}$ is of use. Indeed, our computations show that the corresponding recurrence relation \eqref{dV-Z1-2-Zbar1-2-Zbar2} holds also in this case with $r=1$ and thus, setting ${\rm Re}V_{Z_1^2\oZ_1^2\oZ_2}=0$ enables one to solve it for ${\rm Im}\mu^2_U$. Moreover, the general Lemma \ref{lem-Im-mu2-U} holds also in this branch.

We proceed into the next order with the hope of detecting further general normalizations. Notice that at this stage, the collection of yet unnormalized Maurer-Cartan forms \eqref{MC-3-A''} is reduced to
\begin{equation}\label{MC-5-A''}
 {\rm Re}\mu^1_{Z_1U^l}, \qquad {\rm Im}\mu^2_{Z_1U^l}, \qquad \alpha_{U^{l+1}}
 \end{equation}

In order six, we can still benefit the recurrence relations \eqref{dImV-Z14-Zbar22} with $r=1$ to normalize the Maurer-Cartan forms $\alpha_{UU}$ and ${\rm Im}\mu^1_{Z_1}$ after setting
\[
{\rm Im}V_{Z_1^4\oZ_2^2}={\rm Im}V_{Z_1^3\oZ_1\oZ_2^2}=0.
\]
Indeed, the general Lemma \ref{lem-alpha-UU} holds also in this branch.

Our computations show that on the fiber $\mathcal B_{\bp}$ we have moreover the sixth order relation
\[
dV_{Z_1^2Z_2\oZ_1^3}=-48{\rm i}\,{\rm Im}\mu^2_{Z_1U}+\cdots
\]
where $"\cdots"$ stands for terms which are independent of the Maurer-Cartan form ${\rm Im}\mu^2_{Z_1U}$. Setting
\[
{\rm Im}V_{Z_1^2Z_2\oZ_1^3}=0,
\]
it suggests to solve this recurrence relation for the real Maurer-Cartan form ${\rm Im}\mu^2_{Z_1U}$. More generally we have

\begin{Lemma}
For every $l\geq 0$ and in branch {\rm (A$''$.ii.1)}, one can normalize the Maurer-Cartan form ${\rm Im}\mu^2_{Z_1U^{l+1}}$ by setting ${\rm Im}V_{Z_1^2Z_2\oZ_1^3U^l}=0$.
\end{Lemma}

At this stage, the collection of the remained yet unnormalized Maurer-Cartan forms \eqref{MC-5-A''} is notably reduced to only three real forms
\[
 {\rm Re}\mu^1_{Z_1}, \qquad {\rm Im}\mu^2_{Z_1}, \qquad \alpha_{U},
\]
which, according to \cite{Valiquette-SIGMA}, parameterizes the isotropy group of the real hypersurfaces belonging to this branch. It is in complete agreement with Ershova's result that the isotropy groups associated with these hypersurfaces are of dimensions $\leq 3$. For the {\it model hypersurface}
\[
v=z_1z_2\oz_1+z_1\oz_1\oz_2+z_1^2\oz_2+z_2\oz_1^2,
\]
this group is of the maximum dimension three, generated infinitesimally by the real part of
\begin{equation}\label{field-X}
{\sf X}:={\rm i}\,z_1\partial_{z_2}
\end{equation}
together with the real parts of the two dilation fields \eqref{dilation-fields}.

\begin{Theorem}
Every $2$-nondegenerate real hypersurface $M^5\subset\mathbb C^3$ belonging to the branch {\rm (A$''$.ii.1)} can be transformed to the complete normal form
\begin{equation}\label{NF-A''.ii.1}
 \aligned
 v= z_1z_2\oz_1+z_1\oz_1\oz_2+z_1^2\oz_2+z_2\oz_1^2+V_{Z_2\oZ_2U}\,z_2\oz_2u+\sum_{|\ell_1|+|\ell_2|+l\geq 4}\,\frac{V_{Z^{\ell_1}\oZ^{\ell_2}U^l}}{\ell_1!\,\ell_2!\,l!}\,z^{\ell_1}\oz^{\ell_2} u^l,
 \endaligned
 \end{equation}
 where, regarding the conjugation relation, the coefficients $V_J$ enjoy the cross-section
  \[
 \aligned
 0\equiv V_{Z^\ell U^l}&=V_{Z_1^{j+1}Z_2^{k+1}\oZ_tU^l}=V_{Z_1^{j+3}\oZ_tU^l}={\rm Re}V_{Z_1^2\oZ_1U^l}=V_{Z_1\oZ_tU^{l+1}}
 \\
 &={\rm Re}V_{Z_1^2\oZ_1^2\oZ_2U^l}={\rm Im}V_{Z_1^2\oZ_2U^l}={\rm Im}V_{Z_1^4\oZ_2^2U^l}={\rm Im}V_{Z_1^3\oZ_1\oZ_2^2U^l}={\rm Im}V_{Z_1^2Z_2\oZ_1^3U^l}
 \endaligned
 \]
 for $t=1,2$, $\ell\in\mathbb N_0$ and $i,j,k\geq 0$. Moreover, the isotropy groups of $M^5$ is at most $3$-dimensional.
\end{Theorem}

As in the case of (A$'$.ii.1), we remark that there exists an enormous number of hypersurfaces in this branch with the isotropy groups of dimensions $\lneq 3$. For example, if the third order real invariant $V_{Z_2\oZ_2U}$ does not vanish at $\mathcal B_{\bp}$, then we can set it to $1$ and normalize the Maurer-Cartan form $\alpha_U$. In this case, setting ${\rm Re}V_{Z_1^2Z_2\oZ_1\oZ_2}=0$ enables one to even normalize in addition the Maurer-Cartan form ${\rm Im}\mu^2_{Z_1}$. Accordingly, the isotropy groups of the appearing normal form hypersurfaces
\[
 v= z_1z_2\oz_1+z_1\oz_1\oz_2+z_1^2\oz_2+z_2\oz_1^2+z_2\oz_2u+\sum_{|\ell_1|+|\ell_2|+l\geq 4}\,\frac{V_{Z^{\ell_1}\oZ^{\ell_2}U^l}}{\ell_1!\,\ell_2!\,l!}\,z^{\ell_1}\oz^{\ell_2} u^l
\]
are of dimensions either zero or one. In particular, when the above lifted invariants $V_J$, $\# J\geq 4$ vanish identically, then the isotropy group of the resulted hypersurface
\[
v= z_1z_2\oz_1+z_1\oz_1\oz_2+z_1^2\oz_2+z_2\oz_1^2+z_2\oz_2u
\]
is $1$-dimensional which\,\,---\,\,corresponding to the remained unnormalized real Maurer-Cartan form ${\rm Re}\mu^1_{Z_1}$\,\,---\,\,is generated by the real part of the single infinitesimal generator
\[
{\sf D}_1+2\,{\sf D}_2=z_1\partial_{z_1}+2w\,\partial_w.
\]

\subsection{Branch (A.ii.2)}

As in the branch (A.ii.1), here we have the partially lifted invariant $V_{Z_1Z_2\oZ_1}$ nonzero and we are still permitted to set it to $1$. It results in enjoying as well the normalizations introduced in Lemmas \ref{mu2-Z1-Z2-Zbar1}\,--\,\ref{lem-Im-mu1-Z1-U}. Thus at this point, we may view \eqref{MC-2} as the current collection of the remained unnormalized Maurer-Cartan forms.

We shall assume here that $r=1$ and, in contrary to the subbranch (A$''$.ii.1), the crucial lifted invariant $V_{Z_1^2\oZ_1}$ is now nonzero. By these assumptions, we have the recurrence relation \eqref{dV-z1z1oz1} as
 \[
 dV_{Z_1^2\oZ_1}\equiv V_{Z_1^2\oZ_1}\big(\alpha_U-3\,{\rm Re}\mu^1_{Z_1}\big)-4\,{\rm Re}\mu^2_{Z_1}.
 \]
To simultaneously normalize both Maurer-Cartan forms ${\rm Re}\mu^1_{Z_1}$ and ${\rm Re}\mu^2_{Z_1}$ through the above relation, we set the nonzero lifted invariant $V_{Z_1^2\oZ_1}$ to some imaginary constant, say
\[
V_{Z_1^2\oZ_1}=2{\rm i}.
\]
More generally we have

\begin{Lemma}
In branch {\rm (A.ii.2)} and for each $j,l\geq 0$ with $(j,l)\neq(0,0)$, setting $V_{Z_1^{j+2}\oZ_1U^l}=0$ normalizes identically the two real Maurer-Cartan forms ${\rm Re}\mu^1_{Z_1^{j+1}U^l}$ and ${\rm Re}\mu^2_{Z_1^{j+1}U^l}$.
\end{Lemma}

Based on our computations, it is still possible to normalize ${\rm Re}\mu^2_U$ and $\mu^1_U$ by setting respectively $V_{Z_1\oZ_1U}=0$ and $V_{Z_1\oZ_2U}=0$. More generally, the observations of Lemmas \ref{lem-Re-mu2-U} and \ref{lem-mu1-U} are in turn correct also in this branch.

By the normalizations, applied thus far, the collection \eqref{MC-2} of our remaining unnormalized Maurer-Cartan forms is now reduced to
 \begin{equation}\label{MC-3-A-ii.2}
 {\rm Im}\mu^2_{Z_1^{j}U^l}, \qquad \alpha_{U^{l+1}}, \qquad {\rm for} \ j, l\geq 0.
 \end{equation}

In order to pursue additional normalizations, let us proceed to the next orders.

\subsubsection{Orders four and five of branch {\rm (A.ii.2)}}

In order four and on the fiber $\mathcal B_{\bp}$, we have the recurrence relation
\[
d{\rm Re}V_{Z_1^3\oZ_2}\equiv 3\,{\rm Im}\mu^2_{Z_1^2}+\cdots
\]
where the $"\cdots"$ part is independent of the Maurer-Cartan form ${\rm Im}\mu^2_{Z_1^2}$. Then, by setting
\[
{\rm Re}V_{Z_1^3\oZ_2}=0
\]
it offers to normalize the real form ${\rm Im}\mu^2_{Z_1^2}$. More generally

\begin{Lemma}
For every $j, l\geq 0$ and in branch {\rm (A.ii.2)}, one can normalize the real Maurer-Cartan form ${\rm Im}\mu^2_{Z_1^{j+2}U^l}$ by setting ${\rm Re}V_{Z_1^{j+3}\oZ_2U^l}=0$.
\end{Lemma}

This observation reduces the collection \eqref{MC-3-A-ii.2} of yet unnormalized Maurer-Cartan forms to
\begin{equation}\label{MC-4-A-ii.2}
 {\rm Im}\mu^2_{U^{l+1}}, \qquad {\rm Im}\mu^2_{Z_1U^l}, \qquad \alpha_{U^{l+1}}, \qquad {\rm for} \ l\geq 0.
 \end{equation}

 Similar to the previous branches, we do not find any additional lifted invariants in order four that would allow us to achieve {\it general} normalizations of extra Maurer-Cartan forms. Therefore, we will explore the order five, and our initial candidate for a lifted invariant is ${\rm Re}V_{Z_1^2\oZ_1^2\oZ_2}$, where its recurrence relation on the fiber $\mathcal B_{\bp}$ is
\[
d{\rm Re}V_{Z_1^2\oZ_1^2\oZ_2}\equiv 12\,{\rm Im}\mu^2_U+\cdots.
\]
One can solve this relation for the real Maurer-Cartan form $\mu^2_{U}$ after setting
\[
{\rm Re}V_{Z_1^2\oZ_1^2\oZ_2}=0.
\]
More generally we have

\begin{Lemma}
In branch {\rm (A.ii.2)} and for every $l\geq 0$, one can normalize ${\rm Im}\mu^2_{U^{l+1}}$ and ${\rm Im}\mu^2_{ZU^{l+1}}$ by setting respectively ${\rm Re}V_{Z_1^2\oZ_1^2\oZ_2U^l}=0$ and ${\rm Re}V_{Z_1^3\oZ_1^2\oZ_2U^l}=0$.
\end{Lemma}

Next, we shall look for suitable lifted invariants to normalize the remained real Maurer-Cartan forms ${\rm Im}\mu^2_{Z_1}$ and $\alpha_{U^{l+1}}, l\geq 0$. Unfortunately and as before, we do not find such invariants in order five. Thus, we shall move to order six. We select again in this order the imaginary part of the lifted invariant $V_{Z_1^4\oZ_2^2}$ where its recurrence relation at $\mathcal B_{\bp}$ is
\[
d{\rm Im}V_{Z_1^4\oZ_2^2}\equiv -16{\rm i}\,\alpha_{UU}+\cdots.
\]
Thus, providing ${\rm Im}V_{Z_1^4\oZ_2^2}=0$ enables one to solve the above relation for $\alpha_{UU}$. In general we have

\begin{Lemma}
For every $l\geq 0$ and in branch {\rm (A.ii.2)}, one can normalize the real Maurer-Cartan form $\alpha_{U^{l+2}}$ by setting ${\rm Im}V_{Z_1^4\oZ_2^2U^l}=0$.
\end{Lemma}

By this observation, the two real forms
\[
{\rm Im}\mu^2_{Z_1}, \qquad {\rm and} \qquad \alpha_U
\]
are the only remained Maurer-Cartan forms that did not admit any normalization.
Thus, the dimensions of the isotropy groups associated to the hypersurfaces of branch (A.ii.2) are at most two, verifying Ershova's computations \cite{Ershova-01}. The maximum dimension is enjoyed by the {\it model hypesurface}
\[
v=z_1z_2\oz_1+z_1\oz_1\oz_2+z_1^2\oz_2+z_2\oz_1^2+{\rm i}\,\big(z_1^2\oz_1-z_1\oz_1^2\big)
\]
which admits the real parts of the {\it dilation} and {\it linear} infinitesimal transformations
\begin{equation}\label{infinit-D3-L}
{\sf D}_3:=z_1\partial_{z_1}+z_2\partial_{z_2}+3\,w\partial_{w}, \qquad {\sf L}:={\rm i}\,z_1\partial_{z_2},
\end{equation}
as the generators of its corresponding isotropy algebra.

\begin{Theorem}
Every $2$-nondegenerate real hypersurface $M^5$ of $\mathbb C^3$ belonging to the branch {\rm (A.ii.2)} can be transformed to the complete normal form
\begin{equation}
\label{NF-A-ii.2}
v=z_1z_2\oz_1+z_1\oz_1\oz_2+z_1^2\oz_2+z_2\oz_1^2+{\rm i}\,\big(z_1^2\oz_1-z_1\oz_1^2\big)+V_{Z_2\oZ_2U} z_2\oz_2u+\sum_{|\ell_1|+|\ell_2|+l\geq 4}\,\frac{V_{Z^{\ell_1}\oZ^{\ell_2}U^l}}{\ell_1!\,\ell_2!\,l!}\,z^{\ell_1}\oz^{\ell_2} u^l
\end{equation}
 where, regarding the conjugation relation, the coefficients $V_J$ enjoy the cross-section
 \[
 \aligned
 0\equiv V_{Z^\ell U^l}&=V_{Z_1^{j+1}Z_2^{k+1}\oZ_tU^l}=V_{Z_1^{j+3}Z_2U^l}=V_{Z_1^{j+2}\oZ_1U^l}=V_{Z_1\oZ_tU^{l+1}}
 \\
 &={\rm Im}V_{Z_1^2\oZ_2U^l}={\rm Re}V_{Z_1^2\oZ_1^2\oZ_2U^l}={\rm Re}V_{Z_1^3\oZ_1^2\oZ_2U^l}={\rm Im}V_{Z_1^4\oZ_2^2U^l}
 \endaligned
 \]
 for $t=1,2$, $\ell\in\mathbb N_0$ and $i,j,k\geq 0$. Moreover, the isotropy group of $M^5$ is of dimension $\leq 2$.
\end{Theorem}

As in the former branches, we emphasize that one finds a large number of hypersurfaces in this branch with the isotropy groups of dimensions $\lneqq 2$. For example, if $V_{Z_2\oZ_2U}$  does not vanish on the hypersurface, then setting it to $1$ and ${\rm Re}V_{Z_1^2Z_2\oZ_1\oZ_2}=0$ provides the opportunity of normalizing simultaneously both the remained Maurer-Cartan forms  $\alpha_U$ and ${\rm Im}\mu^2_{Z_1}$. In this case, the resulted normal form hypersurfaces
\[
v=z_1z_2\oz_1+z_1\oz_1\oz_2+z_1^2\oz_2+z_2\oz_1^2+{\rm i}\,\big(z_1^2\oz_1-z_1\oz_1^2\big)+ z_2\oz_2u+\sum_{|\ell_1|+|\ell_2|+l\geq 4}\,\frac{V_{Z^{\ell_1}\oZ^{\ell_2}U^l}}{\ell_1!\,\ell_2!\,l!}\,z^{\ell_1}\oz^{\ell_2} u^l
\]
admit just a trivial isotropy group at $\bp$.

\section{Branch (A.{\rm ii}.3)}
\label{sec-A.ii.3}

In this branch, we assume\,\,---\,\,in contrary to (A.ii.1) and (A.ii.2)\,\,---\,\,that $\Delta_{12}$ vanishes at $\bp$ but
\[
\Delta_{23}=v_{z_1z_2\oz_2}\cdot v_{z_2^2\oz_1}-v_{z_1z_2\oz_1}\cdot v_{z_2^2\oz_2}
\]
is nonzero at this point. Therefore, at least one of the two multiplications $v_{z_1z_2\oz_2}\cdot v_{z_2^2\oz_1}$ and $v_{z_1z_2\oz_1}\cdot v_{z_2^2\oz_2}$ is nonzero at $\bp$. If the former multiplication does not vanish at  this point, then by simply interchanging the role of $z_1$ and $z_2$, we can assume instead that $v_{z_1z_2\oz_1}$ and $v_{z_1^2\oz_2}$ are nonzero. This, brings us back to the former cases (A.ii.1) and (A.ii.2). Then, in this branch, we assume that the jet coordinates $v_{z_1z_2\oz_1}$ and $v_{z_2^2\oz_2}$  (and their lifts) do not vanish at $\bp$.
In addition, here we assume that $v_{z_1^2\oz_2}$ and its lifts are zero at $\bp$ since otherwise, we revert back again to the previous branches (A.ii.1) and (A.ii.2).

Proceeding along the same argument to the proof of Proposition \ref{propo-invariancy-2-nondeg}, one proves also in this case that the two combinations $\Delta_{12}$ and $\Delta_{13}$ vanish identically at $\bp$. Then $\Delta_{23}$ remains nonzero under holomorphic transformations.

We continue the normalizations by considering the recurrence relation of $V_{Z_1Z_2\oZ_1}$ in \eqref{rec-rel-ord-3-first-part}. By the assumptions of this branch, this lift of $v_{z_1z_2\oz_1}$ shall remain nonzero and thus we can set it to $1$ in order to normalize the Maurer-Cartan form $\mu^2_{Z_2}$. Thus, similar to the former branches (A.ii.1) and (A.ii.2), here the normalizations observed in Lemma \ref{mu2-Z1-Z2-Zbar1} are available.

Next, setting $V_{Z_1Z_2\oZ_2}=0$, the corresponding recurrence relation on the fiber $\mathcal B_{\bp}$ is
\[
0=dV_{Z_1Z_2\oZ_2}\equiv -\omu^1_{\oZ_2}-V_{Z_2^2\oZ_2}\,\mu^2_{Z_1}.
\]
We can solve this equation to normalize the Maurer-Cartan form $\omu^1_{\oZ_2}$. More generally

\begin{Lemma}
For every $j,k,l\geq 0$ and in branch {\rm (A.ii.3)}, one can normalize the Maurer-Cartan form $\omu^1_{\oZ_1^j\oZ_2^{k+1}U^l}$ by setting $V_{Z_1Z_2\oZ_1^j\oZ_2^{k+1}U^l}=0$.
\end{Lemma}

Next, let us consider the recurrence relation of $V_{Z_1^2\oZ_1}$. On the fiber $\mathcal B_{\bp}$, we have
\[
dV_{Z_1^2\oZ_1}\equiv V_{Z_1^2\oZ_1}\,\big(\alpha_U-2\,\mu^1_{Z_1}-\omu^1_{\oZ_1}\big)-2\,\mu^2_{Z_1}.
\]
It suggests to set $V_{Z_1^2\oZ_1}=0$ and plainly normalize the Maurer-Cartan form $\mu^2_{Z_1}$. In general we have

\begin{Lemma}
In branch {\rm (A.ii.3)} and for every $j, l\geq 0$, one can normalize the Maurer-Cartan form $\mu^2_{Z_1^{j+1}U^l}$ by setting $V_{Z_1^{j+2}\oZ_1U^l}=0$.
\end{Lemma}

Let us examine the recurrence relation of $V_{Z_2^2\oZ_1}$ which on the fiber $\mathcal B_{\bp}$ is now of the form
\[
dV_{Z_2^2\oZ_1}\equiv V_{Z_2^2\oZ_1}\,\big(2\,\mu^1_{Z_1}+\omu^1_{\oZ_1}-\alpha_U\big).
\]
Thus, after applying the above normalizations, the pseudo-group acts by scaling on $V_{Z_2^2\oZ_1}$, when we restrict it to $\mathcal B_{\bp}$. In order to keep the shape of the expressions in branch (A.ii.3), we let it to be nonzero at $\bp$ and thus we may specify
\[
V_{Z_2^2\oZ_1}=2
\]
which provides the normalization of the Maurer-Cartan form $\mu^1_{Z_1}$ by solving the above recurrence relation. More generally, we have

\begin{Lemma}
For every $j,l\geq 0$ and in branch {\rm (A.ii.3)}, solving the recurrence relation of $V_{Z_1^jZ_2^2\oZ_1U^l}=0$ provides the normalization of the Maurer-Cartan form $\mu^1_{Z_1^{j+1}U^l}$.
\end{Lemma}

Since $V_{Z_2^2\oZ_2}$ is nonzero by assumption, one might consider it as the next candidate for normalizations. But, somewhat surprisingly, our computations indicate that modulo the horizontal coframe we have
\[
dV_{Z_2^2\oZ_2}\equiv 0.
\]
In other words, after applying the above normalizations, the lifted invariant $V_{Z_2^2\oZ_2}$ is now independent of the remaining unnormalized group parameters. Thus, holomorphic transformations are ineffective on its value at $\bp$. As \cite{Ebenfelt-98}, we denote
\[
V_{Z_2^2\oZ_2}(\bp)=\lambda
\]
for some constant $0\neq\lambda\in\mathbb C$, {\it uniquely determined} by $V_{Z_2^2\oZ_2}$.

Similarly, the recurrence relation of $V_{Z_1^2\oZ_2}$ on the fiber $\mathcal B_{\bp}$ is now of the form
\[
dV_{Z_1^2\oZ_2}\equiv 0.
\]
Therefore, $V_{Z_1^2\oZ_2}$ is independent of any group parameter and its value at $\bp$, that we assumed to be zero in this branch, remains invariant.

Proceeding further in order three, let us consider the recurrence relations \eqref{rec-rel-ord-3-compl}. According to our computations and on the fiber $\mathcal B_{\bp}$, they are now of the form
\begin{equation}\label{rec-A.ii.3-completed}
\aligned
dV_{Z_1\oZ_1U}&\equiv -\frac{2}{3}\,V_{Z_1\oZ_1U}\,\alpha_U-2\,{\rm Re}\mu^2_U,
\\
dV_{Z_1\oZ_2U}&\equiv -\frac{2}{3}\,V_{Z_1\oZ_2U}\,\alpha_U-\omu^1_U-2\,\omu^2_U,
\\
dV_{Z_2\oZ_2U}&\equiv -\frac{2}{3}\,V_{Z_2\oZ_2U}\,\alpha_U-\lambda\,\mu^2_U-\overline\lambda\,\omu^2_U.
\endaligned
\end{equation}
It is clear that setting
\[
V_{Z_1\oZ_1U}=V_{Z_1\oZ_2U}=0
\]
enables one to solve the first two relations for the Maurer-Cartan forms ${\rm Re}\mu^2_U$ and $\omu^1_U$. More generally we have

\begin{Lemma}
Let $l\geq 0$. Then in branch {\rm (A.ii.3)},
\begin{itemize}
  \item[$1)$] one can normalize the real Maurer-Cartan form ${\rm Re}\mu^2_{U^{l+1}}$ by setting $V_{Z_1\oZ_1U^{l+1}}=0$.
  \item[$2)$] one can normalize the complex Maurer-Cartan form $\omu^1_{U^{l+1}}$ by setting $V_{Z_1\oZ_2U^{l+1}}=0$.
\end{itemize}
\end{Lemma}

Application of the third recurrence relation in \eqref{rec-A.ii.3-completed} in normalizing the Maurer-Cartan form ${\rm Im}\mu^2_U$ requires that ${\rm Im}\lambda\neq 0$. However, with the aim of generality, we proceed our computations without making any extra assumption regarding the value of $\lambda$.

After applying the above normalizations, the collection of the remained unnormalized Maurer-Cartan forms is now reduced to
\begin{equation}\label{MC-A.ii.3-2}
{\rm Im}\mu^2_{U^{l+1}}, \qquad \alpha_{U^{l+1}}, \qquad {\rm for} \ l\geq 0.
\end{equation}

Unfortunately in order four we will not find appropriate invariants for further normalizations. Then, let us move to order five, where one finds the recurrence relation of ${\rm Re}V_{Z_1^2Z_2\oZ_1^2}$ on the fiber $\mathcal B_{\bp}$ in the simple form
\[
d{\rm Re}V_{Z_1^2Z_2\oZ_1^2}\equiv 8\,{\rm Im}\mu^2_U+\cdots.
\]
This relation readily provides the normalization of the real Maurer-Cartan form ${\rm Im}\mu^2_U$ if we specify ${\rm Re}V_{Z_1^2Z_2\oZ_1^2}=0$. More generally,

\begin{Lemma}
For every $l\geq 0$ and in branch {\rm (A.ii.3)}, the real Maurer-Cartan form ${\rm Im}\mu^2_{U^{l+1}}$ can be normalized by solving the recurrence relation of ${\rm Re}V_{Z_1^2Z_2\oZ_1^2U^l}=0$.
\end{Lemma}

It remains now to normalize the Maurer-Cartan forms $\alpha_{U^{l+1}}$ for $l\geq 0$. For this purpose, we have to proceed to the next order six where, on the fiber $\mathcal B_{\bp}$, we have the recurrence relation of the imaginary part of $V_{Z_1Z_2^3\oZ_1^2}$ simply as
\[
d{\rm Im}V_{Z_1Z_2^3\oZ_1^2}\equiv -8\,\alpha_{UU}+\cdots.
\]
Clearly by setting ${\rm Im}V_{Z_1Z_2^3\oZ_1^2}=0$, one can solve the above relation to normalize the Maurer-Cartan form $\alpha_{UU}$. It leads us to the general

\begin{Lemma}
In branch {\rm (A.ii.3)} and for each $l\geq 0$, one can normalize the real Maurer-Cartan form $\alpha_{U^{l+2}}$ by setting ${\rm Im}V_{Z_1Z_2^3\oZ_1^2U^l}$.
\end{Lemma}

At this stage, the real Maurer-Cartan form $\alpha_U$ is the only remaining unnormalized form. Then, as is realized by Ershova in \cite{Ershova-01}, the dimension of the isotropy group at $\bp$ of every real hypersurface in this branch does not exceed one. When all the remained unconsidered lifted invariants $V_J$ vanish identically at $\bp$, normalization of $\alpha_U$ will be certainly impossible and in this case the isotropy group of the appearing {\it model hypersurfaces}
\[
v=z_1z_2\oz_1+z_1\oz_1\oz_2+z_2^2\oz_1+z_1\oz_2^2+\lambda\,z_2^2\oz_2+\overline\lambda\, z_2\oz_2^2, \qquad 0\neq\lambda\in\mathbb C
\]
has the maximum dimension one, generated infinitesimally by the real part of the single dilation ${\sf D}_3$ in \eqref{infinit-D3-L}.

\begin{Theorem}
Every $5$-dimensional $2$-nondegenerate real hypersurface $M^5\subset\mathbb C^3$ belonging to branch {\rm (A.ii.3)} can be transformed to the complete normal form
\begin{equation*}
v=z_1z_2\oz_1+z_1\oz_1\oz_2+z_2^2\oz_1+z_1\oz_2^2+\lambda\,z_2^2\oz_2+\overline\lambda\, z_2\oz_2^2+V_{Z_2\oZ_2U} z_2\oz_2u+\sum_{|\ell_1|+|\ell_2|+l\geq 4}\,\frac{V_{Z^{\ell_1}\oZ^{\ell_2}U^l}}{\ell_1!\,\ell_2!\,l!}\,z^{\ell_1}\oz^{\ell_2} u^l
\end{equation*}
for a unique nonzero integer $\lambda\in\mathbb C$ where, regarding the conjugation relation, the coefficients $V_J$ enjoy the cross-section
 \[
 \aligned
 0\equiv V_{Z^\ell U^l}&=V_{Z_1^{j+1}Z_2^{k+1}\oZ_1U^l}=V_{Z_1^{j+2}\oZ_1U^l}=V_{Z_1Z_2\oZ_1^j\oZ_2^{k+1}U^l}=V_{Z_1^jZ_2^2\oZ_1U^{l}}
 \\
 &=V_{Z_1\oZ_tU^{l+1}}={\rm Re}V_{Z_1^2Z_2\oZ_1^2U^l}={\rm Im}V_{Z_1Z_2^3\oZ_1^2U^l}
 \endaligned
 \]
 for $\ell\in\mathbb N_0$, $t=1,2$ and $i,j,k\geq 0$. Furthermore, the isotropy group associated to $M^5$ at $\bp$ is at most $1$-dimensional.
\end{Theorem}

As before, we remark that there exists a large number of hypersurfaces in branch (A.ii.3) with the trivial isotropy group at $\bp$. For instance, if $V_{Z_2\oZ_2U}$ does not vanish on $\mathcal B_{\bp}$, then we can normalize the only remained Maurer-Cartan form $\alpha_U$ by solving the third recurrence relation in \eqref{rec-A.ii.3-completed} after setting $V_{Z_2\oZ_2U}=1$. In this case, the appearing normal form hypersurfaces
\[
v=z_1z_2\oz_1+z_1\oz_1\oz_2+z_2^2\oz_1+z_1\oz_2^2+\lambda\,z_2^2\oz_2+\overline\lambda\, z_2\oz_2^2+ z_2\oz_2u+\sum_{|\ell_1|+|\ell_2|+l\geq 4}\,\frac{V_{Z^{\ell_1}\oZ^{\ell_2}U^l}}{\ell_1!\,\ell_2!\,l!}\,z^{\ell_1}\oz^{\ell_2} u^l
\]
admit just the trivial isotropy group.

\section{Branches (A.{\rm ii}.4) and (A.{\rm ii}.5)}
\label{sec-A.ii.4-A.ii.5}

In light of their partially normalized defining equations in \eqref{Ebenfelt-Branches}, both branches (A.ii.4) and (A.ii.5) belong the case in which the two combinations $\Delta_{12}$ and $\Delta_{23}$ in \eqref{2-nondeg-cond-2} vanish at $\bp$ but
\[
\Delta_{13}=v_{z_1^2\oz_1}\cdot v_{z_2^2\oz_2}-v_{z_1^2\oz_2}\cdot v_{z_2^2\oz_1}
\]
is nonzero at this point. Proceeding along the same lines as the proof of Proposition \ref{propo-invariancy-2-nondeg}, one shows that this scenario is invariant under holomorphic transformations.

Clearly, in order to have $\Delta_{13}$ nonzero at $\bp$, at least one of the multiplications $v_{z_1^2\oz_1}\cdot v_{z_2^2\oz_2}$ and $v_{z_1^2\oz_2}\cdot v_{z_2^2\oz_1}$ shall be nonzero at this point. The branch (A.ii.4) concerns the case of which the former multiplication (and its lifts) does not vanish while (A.ii.5) considers the second possibility, assuming that the former multiplication vanishes at $\bp$. We continue this section by studying first the branch (A.ii.4).

\subsection{Branch (A.ii.4)}

In this branch, the opportunity of having the lifts $V_{Z_1^2\oZ_1}$ and $V_{Z_2^2\oZ_2}$ nonzero enables us to set them to some nonzero constant, say
\[
V_{Z_1^2\oZ_1}=V_{Z_2^2\oZ_2}=2.
\]
These specifications enables us to normalize the two Maurer-Cartan forms $\mu^1_{Z_1}$ and $\mu^2_{Z_2}$ by solving the first and last recurrence relations of the list \eqref{rec-rel-ord-3-first-part}, respectively. More generally we have

\begin{Lemma}
In branch {\rm (A.ii.4)} and for each $j,k,l\geq 0$ with $(j,k,l)\neq (0,0,0)$,
\begin{itemize}
  \item[$1)$] one can normalize the Maurer-Cartan form $\omu^1_{\oZ_1^{j+1}\oZ_2^kU^l}$ by setting $V_{Z_1^2\oZ_1^{j+1}\oZ_2^kU^l}=0$.
  \item[$2)$] one can normalize the Maurer-Cartan form $\omu^2_{\oZ_1^j\oZ_2^{k+1}U^l}$ by setting $V_{Z_2^2\oZ_1^{j}\oZ_2^{k+1}U^l}=0$.
\end{itemize}
\end{Lemma}

Now, let us consider the following two recurrence relations on the fiber $\mathcal B_{\bp}$
\[
\aligned
dV_{Z_1Z_2\oZ_1}&\equiv \big(\frac{1}{6}\,(4\,V_{Z_1Z_2\oZ_2}-V_{Z_1\oZ_2^2})\,V_{Z_1Z_2\oZ_1}-2\big)\,\mu^1_{Z_2}+\frac{1}{3}\,(V_{Z_2^2\oZ_1}-V_{Z_2\oZ_1\oZ_2})\,V_{Z_1Z_2\oZ_1}\,\omu^1_{\oZ_2}+\cdots,
\\
dV_{Z_1Z_2\oZ_2}&\equiv \big(\frac{1}{6}\,(4\,V_{Z_1Z_2\oZ_1}-V_{Z_2\oZ_1^2})\,V_{Z_1Z_2\oZ_2}-2\big)\,\mu^2_{Z_1}+\frac{1}{3}\,(V_{Z_1^2\oZ_2}-V_{Z_1\oZ_1\oZ_2})\,V_{Z_1Z_2\oZ_2}\,\omu^2_{\oZ_1}+\cdots,
\endaligned
\]
where $"\cdots"$ represents the terms which do not include the explicitly written Maurer-Cartan forms. Solving these relations results in normalizing the Maurer-Cartan forms $\mu^1_{Z_2}$ and $\mu^2_{Z_1}$, respectively, after setting $V_{Z_1Z_2\oZ_1}=V_{Z_1Z_2\oZ_2}=0$. More generally,

\begin{Lemma}
  Let $j,l\geq 0$. Then, in branch {\rm (A.ii.4)},
  \begin{itemize}
    \item[$1)$] one can normalize the Maurer-Cartan form $\mu^1_{Z_2^{j+1}U^l}$ by setting $V_{Z_1Z_2^{j+1}\oZ_1U^l}=0$.
    \item[$2)$] one can normalize the Maurer-Cartan form $\mu^2_{Z_1^{j+1}U^l}$ by setting $V_{Z_1^{j+1}Z_2\oZ_2U^l}=0$.
  \end{itemize}
\end{Lemma}

Along the way and on the fiber $\mathcal B_{\bp}$, now the two remained recurrence relations in the list \eqref{rec-rel-ord-3-first-part} are simply
\[
dV_{Z_1^2\oZ_2}\equiv 0, \qquad {\rm and} \qquad dV_{Z_2^2\oZ_1}\equiv 0.
\]
Thus, when restricted to the mentioned fiber, these two lifted invariants are now independent of any group parameter. Let us denote
\[
\sigma:=\frac{V_{Z_1^2\oZ_2}}{2}(\bp) \qquad {\rm and} \qquad \nu:=\frac{V_{Z_2^2\oZ_1}}{2}(\bp).
\]
To ensure that the combination $\Delta_{13}$ remains nonzero at point $\bp$, we may assume that $\sigma\cdot\nu\neq 1$.

Now, let us inspect the three remained order three recurrence relations \eqref{rec-rel-ord-3-compl}. On the fiber $\mathcal B_{\bp}$, the first and third relations have now taken the form
\[
dV_{Z_1\oZ_1U}\equiv -4\,{\rm Re}\mu^1_U+\cdots, \qquad dV_{Z_2\oZ_2U}\equiv -4\,{\rm Re}\mu^2_U+\cdots.
\]
Then, by setting $V_{Z_1\oZ_1U}=V_{Z_2\oZ_2U}=0$, one can readily solve these two equations to normalize respectively the real Maurer-Cartan forms ${\rm Re}\mu^1_U$ and ${\rm Re}\mu^2_U$. More generally we have

\begin{Lemma}
For each $l\geq 0$ and in branch {\rm (A.ii.4)},
\begin{itemize}
  \item[$1)$] one can normalize the real Maurer-Cartan form ${\rm Re}\mu^1_{U^{l+1}}$ by setting $V_{Z_1\oZ_1U^{l+1}}=0$.
  \item[$2)$] one can normalize the real Maurer-Cartan form ${\rm Re}\mu^2_{U^{l+1}}$ by setting $V_{Z_2\oZ_2U^{l+1}}=0$.
\end{itemize}
\end{Lemma}

The second recurrence relation in \eqref{rec-rel-ord-3-compl}, namely that of $dV_{Z_1\oZ_2U}$ offers no normalization unless we make the extra assumption that $V_{Z_1\oZ_2U}\neq 0$. Since we do not aim to produce further subbranches, let us search for possible general normalizations in the next orders. Before it, notice that after the above normalizations, the collection of yet unnormalized Maurer-Cartan forms \eqref{MC-1} is now reduced to the real forms
\begin{equation}\label{MC-2-A.ii.4}
{\rm Im}\mu^1_{U^{l+1}}, \qquad {\rm Im}\mu^2_{U^{l+1}}, \qquad \alpha_{U^{l+1}}, \qquad {\rm for} \ l\geq 0.
\end{equation}

According to our computations, at order four, we will find no lifted invariant to provide further normalizations. Therefore, we need to proceed to the next order five where we encounter the following two real recurrence relations on the fiber $\mathcal B_{\bp}$
\begin{equation}\label{rec-ord-5}
\aligned
d{\rm Re}V_{Z_1^2Z_2\oZ_1\oZ_2}&\equiv 4\,|\sigma|^2\,{\rm Im}\mu^1_U+4\,\big(\sigma\nu+2\big)\,{\rm Im}\mu^2_U+\cdots,
\\
d{\rm Re}V_{Z_1Z_2^2\oZ_1\oZ_2}&\equiv 4\,\big(\sigma\nu+2\big)\,{\rm Im}\mu^1_U+4\,|\nu|^2\,{\rm Im}\mu^2_U+\cdots.
\endaligned
\end{equation}
Here the terms in $"\cdots"$ are independent of the Maurer-Cartan forms $\mu^1_U, \mu^2_U$ or their conjugations. By inspecting these relations, one finds out that when $\sigma\nu\neq -1$, then after setting
\[
{\rm Re}V_{Z_1^2Z_2\oZ_1\oZ_2}={\rm Re}V_{Z_1Z_2^2\oZ_1\oZ_2}=0,
\]
the above two recurrence relations provide two linearly independent equations where their solutions give the normalized expressions of ${\rm Im}\mu^1_U$ and $\rm{Im}\mu^2_U$. In general

\begin{Lemma}
Let $l\geq 0$. Then in branch {\rm (A.ii.4)} and with the assumption $\sigma\nu\neq -1$, setting
\[
{\rm Re}V_{Z_1^2Z_2\oZ_1\oZ_2U^l}={\rm Re}V_{Z_1Z_2^2\oZ_1\oZ_2U^l}=0
\]
enables one to normalize the Maurer-Cartan forms ${\rm Im}\mu^1_{U^{l+1}}$ and ${\rm Im}\mu^2_{U^{l+1}}$.
\end{Lemma}

But, when $\sigma\nu=-1$, the above two recurrence relations \eqref{rec-ord-5} will not provide a full rank system for the Maurer-Cartan forms ${\rm Im}\mu^1_U$ and $\rm{Im}\mu^2_U$. In this case, we consider the first recurrence relation together with that of ${\rm Re}V_{Z_1^3\oZ_1\oZ_2}$, which on the fiber $\mathcal B_{\bp}$ give
\[
\aligned
d{\rm Re}V_{Z_1^2Z_2\oZ_1\oZ_2}&\equiv 4\,|\sigma|^2\,{\rm Im}\mu^1_U+4\,{\rm Im}\mu^2_U+\cdots,
\\
d{\rm Re}V_{Z_1^3\oZ_1\oZ_2}&\equiv 36\,\sigma\,{\rm Im}\mu^1_U-\frac{12}{\overline\sigma}\,{\rm Im}\mu^2_U+\cdots.
\endaligned
\]
By setting ${\rm Re}V_{Z_1^2Z_2\oZ_1\oZ_2}={\rm Re}V_{Z_1^3\oZ_1\oZ_2}=0$, one finds the above equations as a full rank homogeneous system which can be solved for  ${\rm Im}\mu^1_U$ and $\rm{Im}\mu^2_U$.

\begin{Lemma}
Let $l\geq 0$. Then in branch {\rm (A.ii.4)} and with the assumption $\sigma\nu= -1$, setting
\[
{\rm Re}V_{Z_1^2Z_2\oZ_1\oZ_2U^l}={\rm Re}V_{Z_1^3\oZ_1\oZ_2U^l}=0
\]
enables one to normalize the Maurer-Cartan forms ${\rm Im}\mu^1_{U^{l+1}}$ and ${\rm Im}\mu^2_{U^{l+1}}$.
\end{Lemma}

To ormalize the remained Maurer-Cartan forms $\alpha_{U^{l+1}}$, $l\geq 0$, we shall move to the next order six. After applying tedious computations, we found the following three recurrence relations on the fiber $\mathcal B_{\bp}$
\begin{equation}\label{ord-6-A.ii.4}
\aligned
dV_{Z_1^2Z_2^2\oZ_1\oZ_2}&\equiv -\frac{8\rm i}{3}\,\big(1+\sigma\nu\big)\,\alpha_{UU}+\cdots,
\\
dV_{Z_1^4\oZ_2^2}&\equiv 16{\rm i}\,\big(\overline\nu-\sigma^2\big)\,\alpha_{UU}+\cdots,
\\
dV_{Z_2^4\oZ_1^2}&\equiv 16{\rm i}\,\big(\overline\sigma-\nu^2\big)\,\alpha_{UU}+\cdots,
\endaligned
\end{equation}
where $"\cdots"$ stands for terms that do not include $\alpha_{UU}$ (or $\alpha_U$). One verifies that the coefficients of the Maurer-Cartan form $\alpha_{UU}$ in these relations can not all be zero simultaneously. Thus, it is always possible to employ one of them for normalizing this form. More precisely, we will have the following four possibilities:

\noindent
$\bullet$ First, if ${\rm Im}(\sigma\nu)\neq 0$, then one can solve the real part of the first recurrence relation of \eqref{ord-6-A.ii.4} for $\alpha_{UU}$ after setting
\[
{\rm Re}V_{Z_1^2Z_2^2\oZ_1\oZ_2}=0.
\]
$\bullet$ Second, if ${\rm Im}(\sigma\nu)=0$ but $\sigma\nu\neq -1$, then one can normalize $\alpha_{UU}$ by solving the imaginary part of the first recurrence relation after setting
\[
{\rm Im}V_{Z_1^2Z_2^2\oZ_1\oZ_2}=0.
\]
$\bullet$ Third, in the case that $\sigma\nu=-1$ but $\overline\nu-\sigma^2\neq 0$, or equivalently when $\sigma\nu=-1$ but $\sigma\neq -1$, then the second recurrence relation in \eqref{ord-6-A.ii.4} offers the normalization of $\alpha_{UU}$ by setting ${\rm Re}V_{Z_1^4\oZ_2^2}=0$ if ${\rm Im}(\overline\nu-\sigma^2)$ is nonzero and by setting ${\rm Im}V_{Z_1^4\oZ_2^2}=0$, otherwise.

\noindent
$\bullet$ Fourth, if $\sigma\nu=-1$ and $\sigma= -1$, or equivalently when $\sigma=-1$ and $\nu=1$, one normalizes $\alpha_{UU}$ by solving the last relation in \eqref{ord-6-A.ii.4} after setting
\[
{\rm Im}V_{Z_2^4\oZ_1^2}=0.
\]

More generally we have that

\begin{Lemma}
Let $l\geq 0$. Then, in branch {\rm (A.ii.4)}, one can normalize the Maurer-Cartan form $\alpha_{U^{l+2}}$ via one of the following ways
\begin{itemize}
  \item[$\rm i)$] by setting ${\rm Re}V_{Z_1^2Z_2^2\oZ_1\oZ_2U^l}=0$ if ${\rm Im}(\sigma\nu)\neq 0$.
  \item[$\rm ii)$] by setting ${\rm Im}V_{Z_1^2Z_2^2\oZ_1\oZ_2U^l}=0$ if ${\rm Im}(\sigma\nu)= 0$ but $\sigma\nu\neq -1$.
  \item[$\rm iii)$] when $\sigma\nu=-1$ but $\sigma\neq -1$, then by setting ${\rm Re}V_{Z_1^4\oZ_2^2U^l}=0$ if ${\rm Im}(\overline\nu-\sigma^2)$ is nonzero and by setting ${\rm Im}V_{Z_1^4\oZ_2^2U^l}=0$, otherwise.
  \item[$\rm iv)$] by setting ${\rm Im}V_{Z_2^4\oZ_1^2U^l}=0$ when $\sigma=-1$ and $\nu=1$.
\end{itemize}
\end{Lemma}

At this stage, except the single real form $\alpha_U$, all the basis Maurer-Cartan forms are normalized in branch (A.ii.4). Thus, the isotropy groups of the real hypersurfaces in this branch are of dimensions at most one. This confirms Ershova's upper bound in \cite{Ershova-01}. The maximum dimension of the isotropy group is enjoyed by the {\it model hypersurfaces}
\[
v=z_1^2\oz_1+z_1\oz_1^2+z_2^2\oz_2+z_2\oz_2^2+\sigma\,z_1^2\oz_2+\overline\sigma\,z_2\oz_1^2+\nu\,z_2^2\oz_1+\overline\nu\,z_1\oz_2^2, \qquad \sigma\nu\neq 1,
\]
which admit the real part of the single dilation ${\sf D}_3$ in \eqref{infinit-D3-L} as the infinitesimal generator of their isotropy algebras at $\bp$.

\begin{Theorem}
Let $M^5$ be a $5$-dimensional $2$-nondegenerate real hypersurface of $\mathbb C^3$ belonging to branch {\rm (A.ii.4)}. Then it can be transformed to the complete normal form
\begin{equation*}
\aligned
v&=z_1^2\oz_1+z_1\oz_1^2+z_2^2\oz_2+z_2\oz_2^2+\sigma\,z_1^2\oz_2+\overline\sigma\,z_2\oz_1^2+\nu\,z_2^2\oz_1+\overline\nu\,z_1\oz_2^2+V_{Z_1\oZ_2U} z_1\oz_2u+V_{Z_2\oZ_1U} z_2\oz_1u
\\
& \ \ \ \ \ \ \ \ \ \ \ \ \ \ \ \ \ \ \ \ \ \ \ \ +\sum_{|\ell_1|+|\ell_2|+l\geq 4}\,\frac{V_{Z^{\ell_1}\oZ^{\ell_2}U^l}}{\ell_1!\,\ell_2!\,l!}\,z^{\ell_1}\oz^{\ell_2} u^l
\endaligned
\end{equation*}
for two unique integers $\sigma, \nu\in\mathbb C$ with $\sigma\nu\neq 1$ where, regarding the conjugation relation, the coefficients $V_J$ enjoy the cross-section
 \[
 \aligned
 0&\equiv V_{Z^\ell U^l}=V_{Z_1^2\oZ_1^{j+1}\oZ_2^kU^l}=V_{Z_2^2\oZ_1^j\oZ_2^{k+1}U^l}=V_{Z_1Z_2^{j+1}\oZ_1U^l}=V_{Z_1^{j+1}Z_2\oZ_2U^{l}}=V_{Z_t\oZ_tU^{l+1}}
 \\
 &={\rm Re}V_{Z_1^2Z_2\oZ_1\oZ_2U^l}={\rm Re}V_{Z_1Z_2^2\oZ_1\oZ_2U^l},
 \endaligned
 \]
 supplemented with
 \[
 \aligned
   \left\{
   \begin{array}{ll}
  {\rm Re}V_{Z_1^2Z_2\oZ_1\oZ_2U^l}=V_{Z_1Z_2^2\oZ_1\oZ_2U^l}=0  & \qquad if \ \ \sigma\nu\neq -1, \\
  {\rm Re}V_{Z_1^2Z_2\oZ_1\oZ_2U^l}={\rm Re}V_{Z_1^3\oZ_1\oZ_2U^l}=0  & \qquad if \ \ \sigma\nu= -1,\\
    \end{array}
   \right.
 \endaligned
 \]
 and
 \[
 \aligned
   \left\{
   \begin{array}{ll}
  {\rm Re}V_{Z_1^2Z_2^2\oZ_1\oZ_2U^l}=0  & \qquad if \ \ {\rm Im}(\sigma\nu)\neq 0, \\
  {\rm Im}V_{Z_1^2Z_2^2\oZ_1\oZ_2U^l}=0  & \qquad if \ \  {\rm Im}(\sigma\nu)=0, \ \sigma\nu\neq -1,\\
  {\rm Re}V_{Z_1^4\oZ_2^2U^l}=0 & \qquad if \ \ \sigma\nu=-1, \ \sigma\neq -1, \ {\rm Im}(\overline\nu-\sigma^2)\neq 0, \\
  {\rm Im}V_{Z_1^4\oZ_2^2U^l}=0 & \qquad if \ \ \sigma\nu=-1, \ \sigma\neq -1, \ {\rm Im}(\overline\nu-\sigma^2)= 0, \\
  {\rm Im}V_{Z_2^4\oZ_1^2U^l}=0 & \qquad  if \ \ \sigma=-1, \ \nu=1,
   \end{array}
   \right.
 \endaligned
 \]
 for $\ell\in\mathbb N_0$, $j,k,l\geq 0$ and $t=1,2$. Furthermore, the isotropy group associated to $M^5$ at $\bp$ is either trivial or $1$-dimensional.
\end{Theorem}

As before, one finds an enormous number of hypersurfaces in this branch admitting trivial isotropy group. For example, if we have ${\rm Re}V_{Z_2\oZ_1U}\neq 0$ then, by setting it to one, we can normalize the only remained Maurer-Cartan form $\alpha_U$ by solving the real part of the third recurrence relation in \eqref{rec-A.ii.3-completed}. Denoting ${\rm Im}V_{Z_2\oZ_1U}(\bp)=\gamma$, then the isotropy groups of the appearing normal form hypersurfaces
\[
\aligned
v&=z_1^2\oz_1+z_1\oz_1^2+z_2^2\oz_2+z_2\oz_2^2+\sigma\,z_1^2\oz_2+\overline\sigma\,z_2\oz_1^2+\nu\,z_2^2\oz_1+\overline\nu\,z_1\oz_2^2
\\
& \ \ \ \ \ \ \ \ \ \ \ \ +(1+{\rm i}\gamma)z_1\oz_2u+(1-{\rm i}\gamma) z_2\oz_1u+\sum_{|\ell_1|+|\ell_2|+l\geq 4}\,\frac{V_{Z^{\ell_1}\oZ^{\ell_2}U^l}}{\ell_1!\,\ell_2!\,l!}\,z^{\ell_1}\oz^{\ell_2} u^l
\endaligned
\]
are nothing but the trivial.

\subsection{Branch (A.ii.5)}

Now we consider the last branch (A.ii.5) in the list \eqref{Ebenfelt-Branches}. As mentioned at the beginning of this section, here we assume that the multiplication of the two partially lifted invariants $V_{Z_1^2\oZ_1}$ and $V_{Z_2^2\oZ_2}$ vanishes, contrary to the multiplication of $V_{Z_1^2\oZ_2}$ and $V_{Z_2^2\oZ_1}$.
We begin the normalizations by setting
\[
V_{Z_1^2\oZ_2}=V_{Z_2^2\oZ_1}=2
\]
as we are permitted here. Then, solving the recurrence relations of these two lifted invariants in \eqref{rec-rel-ord-3-first-part} provides respectively the normalizations of  $\mu^1_{Z_1}$ and $\mu^2_{Z_2}$. In general, we have

\begin{Lemma}
Let $j,k,l\geq 0$ with $(j,k,l)\neq (0,0,0)$. Then in branch {\rm (A.ii.5)},
\begin{itemize}
  \item[$1)$] one can normalize the Maurer-Cartan form $\mu^1_{Z_1^{j+1}Z_2^kU^l}$ by setting $V_{Z_1^{j+2}Z_2^k\oZ_2U^l}=0$.
  \item[$2)$] one can normalize the Maurer-Cartan form $\mu^2_{Z_1^jZ_2^{k+1}U^l}$ by setting $V_{Z_1^jZ_2^{k+2}\oZ_1U^l}=0$.
\end{itemize}
\end{Lemma}

Next, let us consider the recurrence relations of $V_{Z_1Z_2\oZ_1}$ and $V_{Z_1Z_2\oZ_2}$ which now have the following forms on the fiber $\mathcal B_{\bp}$
\[
\aligned
dV_{Z_1Z_2\oZ_1}&\equiv {  A}\,V_{Z_1Z_2\oZ_1}-V_{Z_1^2\oZ_1}\,\mu^1_{Z_2}-V_{Z_1Z_2\oZ_2}\,\omu^2_{\oZ_1}-2\,\mu^2_{Z_1},
\\
dV_{Z_1Z_2\oZ_2}&\equiv { B}\,V_{Z_1Z_2\oZ_2}-V_{Z_2^2\oZ_2}\,\mu^2_{Z_1}-V_{Z_1Z_2\oZ_1}\,\omu^1_{\oZ_2}-2\,\mu^1_{Z_2}.
\endaligned
\]
Here, $A$ and $B$ are two certain polynomials in terms of the Maurer-Cartan forms and lifted invariants. Thus, clearly, one can normalize the Maurer-Cartan forms $\mu^2_{Z_1}$ and $\mu^1_{Z_2}$ by solving respectively the above recurrence relations after setting $V_{Z_1Z_2\oZ_1}=V_{Z_1Z_2\oZ_2}=0$. More generally

\begin{Lemma}
Let $j,l\geq 0$. In branch {\rm (A.ii.5)},
\begin{itemize}
  \item[$1)$] one can normalize the Maurer-Cartan form $\mu^2_{Z_1^{j+1}U^l}$ by setting $V_{Z_1^{j+1}Z_2\oZ_1U^l}=0$.
  \item[$2)$] one can normalize the Maurer-Cartan form $\mu^1_{Z_2^{j+1}U^l}$ by setting $V_{Z_1Z_2^{j+1}\oZ_2U^l}=0$.
\end{itemize}
\end{Lemma}

Along the way, our computations show that on the fiber $\mathcal B_{\bp}$, we simply have
\[
dV_{Z_1^2\oZ_1}\equiv 0, \qquad {\rm and} \qquad dV_{Z_2^2\oZ_2}\equiv 0.
\]
This alludes that the two lifted invariants $V_{Z_1^2\oZ_1}$ and $V_{Z_2^2\oZ_2}$ are now independent of the group parameters. On the other hand, by the assumptions of this branch, we know that the multiplication of $V_{Z_1^2\oZ_1}$ and $V_{Z_2^2\oZ_2}$ vanishes at $\bp$, thus at least one of them is zero at this point. By exchanging the roles of $z_1$ and $z_2$, if necessary, we can always assume that $V_{Z_2^2\oZ_2}(\bp)= 0$. Denote
\[
\eta:=\frac{V_{Z_1^2\oZ_1}(\bp)}{2}.
\]

It remains in order three to consider the recurrence relations \eqref{rec-rel-ord-3-compl}. On the fiber $\mathcal B_{\bp}$, the second relation is now of the form
\[
dV_{Z_1\oZ_2U}\equiv {C}\,V_{Z_1\oZ_2U}-2\,\mu^1_U-2\,\omu^2_U
\]
for some polynomial $C$ in terms of the Maurer-Cartan forms and differential invariants. By setting $V_{Z_1\oZ_2U}=0$, one solves this equation to normalize the Maurer-Cartan form $\mu^1_U$. In general we have

\begin{Lemma}
For every $l\geq 0$ and in branch {\rm (A.ii.5)}, one can normalize the Maurer-Cartan form $\mu^1_{U^{l+1}}$ by solving the recurrence relation of $V_{Z_1\oZ_2U^{l+1}}=0$.
\end{Lemma}

Unfortunately, the other two recurrence relations of $V_{Z_1\oZ_1U}$ and $V_{Z_2\oZ_2U}$ in \eqref{rec-rel-ord-3-compl} are generally not useful for further normalizing the Maurer-Cartan forms. Then, we proceed into the next order four. Before it, we notice that the collection of yet unnormalized basis Maurer-Cartan forms \eqref{MC-1} is now reduced to
\begin{equation}\label{MC-2-A.ii.5}
\mu^2_{U^{l+1}}, \qquad \omu^2_{U^{l+1}}, \qquad \alpha_{U^{l+1}}, \qquad {\rm for} \ l\geq 0.
\end{equation}

In order four, none of the lifted invariants is able to normalize any of these forms. But, in order five and on the fiber $\mathcal B_{\bp}$, we have the  recurrence relation
\[
dV_{Z_2^3\oZ_1^2}\equiv -60{\rm i}\,\mu^2_U+\cdots
\]
where $"\cdots"$ does not comprise $\mu^2_U$ or its conjugation. Then, clearly, setting $V_{Z_2^3\oZ_1^2}=0$ provides the normalization of the Maurer-Cartan form $\mu^2_{U}$. More generally we have

\begin{Lemma}
For $l\geq 0$ and in branch {\rm (A.ii.5)}, one normalizes the Maurer-Cartan form $\mu^2_{U^{l+1}}$ by setting $V_{Z_2^3\oZ_1^2U^l}=0$.
\end{Lemma}

It remains only the normalizations of the real forms $\alpha_{U^{l+1}}$, $l\geq 0$. For this purpose, we have to move to order six, where we have the following recurrence relation on $\mathcal B_{\bp}$

\[
dV_{Z_1^2Z_2^2\oZ_1\oZ_2}\equiv-\frac{8\rm i}{3}\,\alpha_{UU}+\cdots.
\]
Here, $"\cdots"$ stands for the terms which do not admit the Maurer-Cartan form $\alpha_{UU}$ (or $\alpha_U$). Thus, we may normalize this real Maurer-Cartan form by solving the imaginary part of the above relation after setting ${\rm Im}V_{Z_1^2Z_2^2\oZ_1\oZ_2}=0$. More generally, we have

\begin{Lemma}
For every $l\geq 0$ and in branch {\rm (A.ii.5)}, one can normalize the real Maurer-Cartan form $\alpha_{U^{l+2}}$ by setting ${\rm Im}V_{Z_1^2Z_2^2\oZ_1\oZ_2U^l}=0$.
\end{Lemma}

At this stage, $\alpha_U$ is the only remaining unnormalized Maurer-Cartan form in branch (A.ii.5). Therefore, the isotropy groups associated to the hypersurfaces in this branch are of dimensions $\leq 1$. This, confirms Ershova's upper bound in \cite{Ershova-01}. Among the mentioned hypersurfaces, the isotropy groups of the {\it model hypersurfaces}
\[
v=\eta\,z_1^2\oz_1+\overline\eta z_1\oz_1^2+z_1^2\oz_2+z_2\oz_1^2+z_2^2\oz_1+z_1\oz_2^2
\]
are of the maximum dimension one, generated by the real part of the infinitesimal dilation ${\sf D}_3$ in \eqref{infinit-D3-L}.

\begin{Theorem}
Given a $5$-dimensional $2$-nondegenerate real hypersurface $M^5\subset\mathbb C^3$ belonging to the general branch {\rm (A.ii.5)}, there exists an origin-preserving transformation mapping $M^5$ to the complete normal form
\begin{equation*}
\aligned
v&=\eta\,z_1^2\oz_1+\overline\eta z_1\oz_1^2+z_1^2\oz_2+z_2\oz_1^2+z_2^2\oz_1+z_1\oz_2^2+V_{Z_1\oZ_1U} z_1\oz_1u+V_{Z_2\oZ_2U} z_2\oz_2 u
\\
& \ \ \ \ \ \ \ \ \ \ \ \ \ \ \ \ \ \ \ \ \ \ \ \ \ \ \ \ \ +\sum_{|\ell_1|+|\ell_2|+l\geq 4}\,\frac{V_{Z^{\ell_1}\oZ^{\ell_2}U^l}}{\ell_1!\,\ell_2!\,l!}\,z^{\ell_1}\oz^{\ell_2} u^l
\endaligned
\end{equation*}
for some unique integer $\eta\in\mathbb C$ where, regarding the conjugation relation, the coefficients $V_J$ enjoy the cross-section
\[
\aligned
0&\equiv V_{Z^\ell U^l}=V_{Z_1^{j+2}Z_2^k\oZ_2U^l}=V_{Z_1^jZ_2^{k+2}\oZ_1U^l}=V_{Z_1^{j+1}Z_2\oZ_1U^l}=V_{Z_1Z_2^{j+1}\oZ_2U^l}
\\
&=V_{Z_1\oZ_2U^{l+1}}=V_{Z_2^3\oZ_1^2U^l}={\rm Im}V_{Z_1^2Z_2^2\oZ_1\oZ_2U^l}
\endaligned
\]
for $j,k,l\geq 0$. Moreover, the isotropy group of $M^5$ at the point $\bp=0$ is of dimension at most one.
\end{Theorem}

Remark that in branch (A.ii.5), one finds an abundance of hypersurface with trivial isotropy groups. For instance, if either of the order three real invariants ${\bf V}=V_{Z_1\oZ_1U}, V_{Z_2\oZ_2U}$ is nonzero, it becomes possible to normalize $\alpha_U$ solving its corresponding recurrence relation in \eqref{rec-rel-ord-3-compl} after setting ${\bf V}=1$. In this case, the corresponding normal forms
\[
v=\eta\,z_1^2\oz_1+\overline\eta z_1\oz_1^2+z_1^2\oz_2+z_2\oz_1^2+z_2^2\oz_1+z_1\oz_2^2+ z_1\oz_1u+\gamma\, z_2\oz_2 u+\sum_{|\ell_1|+|\ell_2|+l\geq 4}\,\frac{V_{Z^{\ell_1}\oZ^{\ell_2}U^l}}{\ell_1!\,\ell_2!\,l!}\,z^{\ell_1}\oz^{\ell_2} u^l
\]
and
\[
v=\eta\,z_1^2\oz_1+\overline\eta z_1\oz_1^2+z_1^2\oz_2+z_2\oz_1^2+z_2^2\oz_1+z_1\oz_2^2+\gamma\, z_1\oz_1u+ z_2\oz_2 u+\sum_{|\ell_1|+|\ell_2|+l\geq 4}\,\frac{V_{Z^{\ell_1}\oZ^{\ell_2}U^l}}{\ell_1!\,\ell_2!\,l!}\,z^{\ell_1}\oz^{\ell_2} u^l
\]
admit trivial isotropy groups. In these two expressions, $\gamma\in\mathbb R$ is respectively the value of $V_{Z_2\oZ_2U}$ and $V_{Z_1\oZ_1U}$ at the point $\bp$.

\section{The equivalence problem}
\label{sec-equiv-prob}

Now that the complete normal forms associated with all branches \eqref{Ebenfelt-Branches} of the class of $5$-dimensional $2$-nondegenerate hypersurfaces of $\mathbb C^3$ are constructed at Levi non-uniformly rank zero points, we can treat the underlying biholomorphic equivalence. By definition (\cite{Olver-1995}), two submanifolds $M, M'\subset\mathbb C^N$ passing respectively through the distinguished points $p$ and $p'$ are called (locally) {\it biholomorphically equivalent} if there exists some local biholomorphism $\varphi:(\mathbb C^N,p)\rightarrow(\mathbb C^N, p')$ mapping $p$ to $p'$ with $\varphi(M\cap U)=M'\cap U'$ for some local neighborhoods $U, U'\subset\mathbb C^N$ of $p, p'$.

If two real hypersurfaces are biholomorphically equivalent, then they clearly admit the same normal forms. In particular, if they belong to different branches (A.ii.1)--(A.ii.5), they are certainly inequivalent. However, the converse of the mentioned statement is not correct in general. Indeed, if two hypersurfaces admit the same normal forms, then they will be {\it formally} equivalent for sure. Nevertheless, it does not assure their biholomorphic equivalence, in general. Fortunately, in the current case that our hypersurfaces are finitely nondegenerate, we have the following key result

\begin{Theorem}
(\cite[Theorem 5]{BER-97})
Let $M$ and $M'$ be two real-analytic hypersurfaces in $\mathbb C^n$ which both pass through the origin. Assume that $M$ is finitely nondegenerate. Then every origin-preserving formal equivalence map between $M$ and $M'$ is biholomorphic.
\end{Theorem}

Recall that as mentioned in Remark \ref{rem-orig-pres}, our selected zeroth order cross-section imposed the associated normal form transformations being origin-preserving. Accordingly, a straightforward application of the above theorem implies that

\begin{Proposition}
Two hypersurfaces belonging to each of the branches {\rm (A.ii.1)--(A.ii.5)} are equivalent through an origin-preserving  biholomorphism if and only if they admit the same normal form.
\end{Proposition}

Moreover, as a direct consequence of \cite[Theorem 4.11]{Valiquette-SIGMA},  we have

\begin{Proposition}
Let $M$ and $M'$ be two real hypersurfaces in $\mathbb C^3$ belonging to one of the branches  {\rm (A.ii.1)--(A.ii.5)}. Assume that they are biholomorphically equivalent through an origin-preserving holomorphic map $\varphi:M\rightarrow M'$. If $\psi$ is another origin-preserving biholomorphism  between $M$ and $M'$, then there exists two transformations $h$ and $h'$, respectively in the isotropy groups of $M$ and $M'$, enjoying
\[
\psi=h'\circ\varphi\circ h.
\]
\end{Proposition}

\subsection*{Acknowledgments}
%%%%%

The research of the author was supported in part by a grant from IPM, No.\ 14020417.

\vglue10pt
%%%%%

\vglue 10pt


\begin{thebibliography}{99}
%%%%%
\vglue10pt

\bibitem{BER-97}
M.S.\ Baouendi, P.\ Ebenfelt and L.P. Rothschild, Parametrization of local biholomorphisms of real analytic hypersurfaces, {\it Asian J. Math.}, {\bf 1}(1) (1997), 1--16.

\bibitem{BER}
\textrm{M.S. Baouendi, P. Ebenfelt and L.P. Rothschild}, \emph{Real Submanifolds in Complex Space and their Mappings}, Princeton Univ. Press, 1999, xii+404 pp.

\bibitem{Cartan-1932}
\'E. \ Cartan, Sur la g\'eom\'etrie pseudo-conforme des hypersurfaces de l'espace de deux variables complexes, I., {\it Ann. Mat. Pura Appl.}, {\bf 11}(4) (1932), 17--90 and II. {\it Ann. Scuola Norm. Sup. Pisa}, {\bf 1} (2) (1932) 333--354.

\bibitem{Cartan-1935}
\'E. \ Cartan, La M\'ethode du Rep\`{e}re Mobile, la Th\'{e}orie des Groupes Continus, et les Espaces G\'{e}n\'{e}ralis\'{e}s, {\it Expos\'{e}s de G\'{e}om\'{e}trie}, no. 5, Hermann, Paris, 1935.

\bibitem{Chern-Moser}
S.S.\ Chern and J.K.\ Moser, Real hypersurfaces in complex spaces, {\it Acta Math.}, {\bf 133} (1974), 219--271.

\bibitem{Ebenfelt-98}
P. \ Ebenfelt, Normal forms and biholomorphic equivalence of real hypersurfaces in $\mathbb C^3$, {\it Indiana Univ. Math. J.}, {\bf 47}(2) (1998), 311--366.

\bibitem{Ebenfelt-01}
P. \ Ebenfelt, Uniformly Levi degenerate CR manifolds: the $5$-dimensional case, {\it Duke Math. J.}, {\bf 110}(1) (2001), 37--80.

\bibitem{Ershova-01}
A.E. \ Ershova, Automorphisms of $2$-nondegenerate hypersurfaces in $\mathbb C^3$, {\it Math. Notes}, {\bf 69}(2) (2001), 188--195.

\bibitem{Fels-Kaup-08}
G.\ Fels and W.\ Kaup, Classification of Levi degenerate homogeneous CR-manifolds in dimension $5$, {\it Acta Math.}, {\bf 201}(1) (2008), 1--82.

\bibitem{Olver-Fels-97}
M.\ Fels and P.J.\ Olver, On relative invariants, {\it Math. Annal.}, {\bf 308} (1997), 701--732.

\bibitem{Olver-Fels-99}
M.\ Fels and P.J.\ Olver, Moving coframes: II. regularization and theoretical foundations, {\it Acta Appl.\ Math.}, {\bf 55} (1999), 127--208.

\bibitem{FMT-22}
W.G. \ Foo, J. \ Merker and T.A.\ Ta, On convergent Poincar\'{e}-Moser reduction for Levi degenerate embedded $5$-dimensional CR manifolds, {\it New York J. Math.}, {\bf 28} (2022), 250--336.

\bibitem{Isaev-Zaitsev-13}
A.\ Isaev and D.\ Zaitsev, Reduction of five-dimensional uniformly degenerate Levi CR structures to absolute parallelisms, {\it J. Geom. Anal.}, {\bf 23} (2013), 1571--1605.

\bibitem{KK-22}
M.\ Kol\'{a}\v{r} and I. Kossovskiy, A complete normal form for everywhere Levi-degenerate hypersurfaces in $\mathbb C^3$, {\it Adv.\ Math.}, {\bf 408} (2022), 108590.

\bibitem{KKZ-17}
M.\ Kol\'{a}\v{r}, I. Kossovskiy, and D. Zaitsev, Normal forms in Cauchy-Riemann geometry, {\it Contem.\ Math.}, {\bf 681} (2017), 153--177.

\bibitem{Medori-Spiro-14}
C.\ Medori and A.\ Spiro, The equivalence problem for $5$-dimensional Levi degenerate CR manifolds, {\it Int. Math. Res. Not. IMRN}, {\bf 20} (2014), 5602--5647.

\bibitem{Merker-Pocchiola-20}
J.\ Merker and S.\ Pocchiola, Explicit absolute parallelism for 2-nondegenerate real hypersurfaces $M^5\subset\mathbb C^3$ of constant Levi rank $1$, {\it J. Geom. Anal.}, {\bf 30} (2020), 2689--2730.

\bibitem{MPS-13}
J.\ Merker, S.\ Pocchiola and M. \ Sabzevari, Equivalences of 5-dimensional CR-manifolds
II: General classes {\sf I}, {\sf II}, ${\sf III}_1$, ${\sf III}_2$, ${\sf IV}_1$, ${\sf IV}_2$, available at {\tt arxiv.org}, 2013.

\bibitem{Olver-1995}
P.J.\ Olver, \emph{Equivalence, Invariants and Symmetry}, Cambridge, Cambridge University Press, 1995, xvi+525~pp.

\bibitem{Olver-2015}
P.J.\ Olver, Modern developments in the theory and applications of moving frames, {\it London Math.\ Soc. Impact150 Stories} {\bf 1} (2015), 14--50.

\bibitem{Olver-2018}
P.J.\ Olver, Normal forms for submanifolds under group actions, {\it Sym.\ Diff.\ Eq.\ Appl.} (V. G.\ Kac et al.\ (eds.)), Proc.\ Math.\ Stat., Springer (2018), 3--27.

\bibitem{Olver-Pohjanpelto-05}
P.J.\ Olver and J.\ Pohjanpelto, Maurer--Cartan forms and the structure of Lie pseudo groups, {\it Sel.\ Math., New Ser.}, {\bf 11} (2005), 99--126.

\bibitem{Olver-Pohjanpelto-08}
P.J.\ Olver and J.\ Pohjanpelto, Moving frames for Lie pseudo-groups, {\it Canad.\ J.\ Math.}, {\bf 60}(2) (2008), 1336--1386.

\bibitem{OSV-23}
P.J.\ Olver, M.\ Sabzevari and F. \ Valiquette, Normal forms, moving frames and differential invariants for nondegenerate hypersurfaces in $\mathbb C^2$, {\it J. Geom. Anal.}, {\bf 33} (2023), 192.

\bibitem{Poincare-1907}
H.\ Poincar\'e, Les fonctions analytiques de deux variables et la repr\'esentation conforme, {\it Rend.\ Circ.\ Mat.\ Palermo} {\bf 23} (1907), 185--220.

\bibitem{Normal-Form}
M.\ Sabzevari, Convergent normal forms for five dimensional totally nondegenerate CR manifolds in $\mathbb C^4$, {\it J.\ Geom.\ Anal.}, {\bf 31} (2021), 7900--7946.

\bibitem{Maple}
M.\ Sabzevari, Maple worksheet of the corresponding computations, available on demand.

\bibitem{Tanaka-76}
N.\ Tanaka, On non-degenerate real hypersurfaces, graded Lie algebras and Cartan connections, {\it Japan.\ J. \ Math.} {\bf 2} (1976), 131--190.

\bibitem{Valiquette-SIGMA}
F.\ Valiquette, Solving local equivalence problems with the equivariant moving frame method, {\it Sym.\ Int.\ Geom.: Meth.\ Appl.} ({\it SIGMA}), {\bf 9} 029 (2013), 43 pp.



\end{thebibliography}
\end{document}